\documentclass{amsart}
\usepackage{fullpage}
\usepackage{color}
\usepackage{amsmath}
\usepackage{amsthm}
\usepackage{amsfonts}
\usepackage{amscd}

\newtheorem{theorem}{Theorem}[section]
\newtheorem{lemma}[theorem]{Lemma}
\newtheorem{proposition}[theorem]{Proposition}
\newtheorem{corollary}[theorem]{Corollary}
\newtheorem{definition}[theorem]{Definition}

\theoremstyle{remark}
\newtheorem*{remark}{Remark}

\newcounter{puntero}
\newcommand{\Prim}{\mathop{\mathrm{Prim}}}
\newcommand{\Hom}{\mathop{\mathrm{Hom}}}

\newcommand{\mR}{\mathcal{R}}
\newcommand{\mRhat}{\widehat{\mathcal{R}}}

\newcommand{\Id}{\mathrm{Id}}
\newcommand{\bx}{\boldsymbol{x}}
\newcommand{\by}{\boldsymbol{y}}
\newcommand{\bz}{\boldsymbol{z}}
\newcommand{\bzero}{\mathbf{0}}
\newcommand{\be}{\mathbf{e}}

\font\cyrillic=wncyi10
\newcommand{\SU}{\mathop{\hbox{{\cyrillic UX}}}}

\begin{document}

\title[Formal multiplications and non-associative Lie theory]{Formal
multiplications, bialgebras of distributions and non-associative Lie
theory}
\author{J.\ Mostovoy}
\address{Departamento de Matem\'aticas, CINVESTAV, Apartado Postal 14-740
C.P.\ 07000, M\'exico, D.F., M\'exico}
\email{jacob@math.cinvestav.mx}
\author{J.\ M.\ P\'erez--Izquierdo}
\address{Departamento Matem\'aticas y Computaci\'on, Universidad de La
Rioja, Logro\~no, 26004, Spain}
\email{jm.perez@dmc.unirioja.es}
\date{}

\thanks{2000 Mathematics Subject Classification: primary 17D99, secondary 22E60, 20N05. }

\begin{abstract}
We describe the general non-associative version of Lie theory that
relates unital formal multiplications (formal loops), Sabinin
algebras and non-associative bialgebras.

Starting with a formal multiplication we construct a non-associative
bialgebra, namely, the bialgebra of distributions with the
convolution product. Considering the primitive elements in this
bialgebra gives a functor from formal loops to Sabinin algebras. We
compare this functor to that of Mikheev and Sabinin and show that
although the brackets given by both constructions coincide, the
multioperator does not. We also show how identities in loops produce
identities in bialgebras. While associativity in loops translates
into associativity in algebras, other loop identities (such as the
Moufang identity) produce new algebra identities. Finally, we define
a class of unital formal multiplications for which Ado's theorem
holds and give examples of formal loops outside this class.

A by-product of the constructions of this paper is a new identity on
Bernoulli numbers. We give two proofs: one coming from the formula
for the non-associative logarithm, and the other (due to D.\ Zagier)
using generating functions.

\end{abstract}

\maketitle

%\tableofcontents

%%%%%%%%%%%%%%%%%%%%%%%%%%%%%%%%%%%%%%%%%%%%%%%%%%%%%%%%%%%%%%%%%%%%%%%%%%%%%%%%
%%%%%%%%%%%%%%%%%%%%%%%%%%%%%%%%%%%%%%%%%%%%%%%%%%%%%%%%%%%%%%%%%%%%%%%%%%%%%%%%
%%%%%%%%%%%%%%%%                                                %%%%%%%%%%%%%%%%
%%%%%%%%%%%%%%%%                                                %%%%%%%%%%%%%%%%
%%%%%%%%%%%%%%%%                   INTRODUCTION                 %%%%%%%%%%%%%%%%
%%%%%%%%%%%%%%%%                                                %%%%%%%%%%%%%%%%
%%%%%%%%%%%%%%%%                                                %%%%%%%%%%%%%%%%
%%%%%%%%%%%%%%%%%%%%%%%%%%%%%%%%%%%%%%%%%%%%%%%%%%%%%%%%%%%%%%%%%%%%%%%%%%%%%%%%
%%%%%%%%%%%%%%%%%%%%%%%%%%%%%%%%%%%%%%%%%%%%%%%%%%%%%%%%%%%%%%%%%%%%%%%%%%%%%%%%

\section{Introduction}

The Lie theory describes the relationship among three types of
algebraic structures: Lie groups, Lie algebras and Hopf algebras. In
brief, we have the following triangle: for a finite-dimensional Lie
group $G$, the Hopf algebra of distributions on $G$ supported at the
unit is nothing else but the universal enveloping algebra of the Lie
algebra of $G$.

Strictly speaking, Lie algebras correspond directly not to Lie
groups, but rather to {\em formal} Lie groups (for instance, via the
Campbell-Baker-Hausdorff formula). The local methods of Lie theory
are not sufficient to establish that finite-dimensional formal
groups give rise to Lie groups: this follows from the existence of a
faithful representation for every finite-dimensional Lie algebra.
This phenomenon is even more apparent in the Lie theory of
non-associative multiplications, where many natural examples of
multiplications on manifolds are of local nature and do not have
evident extensions to global operations.

In the present paper we describe the non-associative version of the
correspondence between formal groups, Lie algebras and Hopf
algebras.

The main step towards generalizing the Lie theory to this context
was done by Sabinin and Mikheev \cite{SM87, MS90} who defined
algebraic structures tangent to arbitrary local analytic loops
(multiplications). These structures, now known as Sabinin algebras,
can be integrated under some convergence conditions to local loops:
essentially, they are the analog of Lie algebras in the
non-associative setting. Shestakov and Umirbaev later showed
\cite{SU02} that the set of primitive elements in any bialgebra has
the structure of a Sabinin algebra, and it was proved by the second
author of the present paper \cite{Pe07}, that each Sabinin algebra
arises in this way. The main purpose of the present paper is to show
how the Lie theory for non-associative formal multiplications can be
constructed by first passing from a formal multiplication to the
corresponding bialgebra of distributions, and then to the Sabinin
algebra of the primitive elements of the latter. We compare this
construction to the direct geometric argument of Sabinin and Mikheev
and show that these two constructions do not give precisely the same
result: they produce Sabinin algebras with coinciding brackets but
different multioperators.

There are two aspects of the non-associative Lie theory that are
absent from the usual Lie theory. Firstly, in the non-associative
context it is rather usual to consider a class of multiplications
satisfying certain identity. We show how these identities translate
into identities in the corresponding bialgebras of distributions.
The second novelty is that while Lie groups are always locally
isomorphic to linear groups, this property (or, rather, its
appropriate generalization) no longer holds for general loops. We
discuss this phenomenon and give examples of formal loops that do
not satisfy this property.

%All constructions of the present paper are carried out over a field
%of characteristic 0. While it may not be difficult to describe the
%correspondence between formal multiplications and bialgebras in
%arbitrary characteristic, at the moment there is no description of
%the primitive elements functor if the characteristic of the ground
%field is positive.

The paper has the following structure. In the next section we show
that the category of unital formal multiplications is equivalent to
that of irreducible cocommutative and coassociative bialgebras. In
Section~3 we consider how identities in formal loops correspond to
identities in bialgebras. In Section~4 we show that the primitive
operations of Shestakov and Umirbaev give an equivalence between the
category of irreducible cocommutative and coassociative bialgebras
and the category of Sabinin algebras. Section~5 contains a
comparison of two functors from formal loops to Sabinin algebras:
the Sabinin algebra of the primitive elements in the algebra of
distributions on a formal loop and the the Sabinin algebra as
defined by Sabinin and Mikheev. In Section~6 we discuss linear
formal loops (those for which Ado's theorem holds). Finally, in the
appendix we give the formulae for the non-associative exponential
and logarithm and describe an identity on Bernoulli numbers.

%%%%%%%%%%%%%%%%%%%%%%%%%%%%%%%%%%%%%%%%%%%%%%%%%%%%%%%%%%%%%%%%%%%%%%%%%%%%%%%%
%%%%%%%%%%%%%%%%%%%%%%%%%%%%%%%%%%%%%%%%%%%%%%%%%%%%%%%%%%%%%%%%%%%%%%%%%%%%%%%%
%%%%%%%%%%%%%%%%                                                %%%%%%%%%%%%%%%%
%%%%%%%%%%%%%%%%                                                %%%%%%%%%%%%%%%%
%%%%%%%%%%%%%%%%     FORMAL MULTIPLICATIONS AND BIALGEBRAS      %%%%%%%%%%%%%%%%
%%%%%%%%%%%%%%%%                OF DISTRIBUTIONS                %%%%%%%%%%%%%%%%
%%%%%%%%%%%%%%%%                                                %%%%%%%%%%%%%%%%
%%%%%%%%%%%%%%%%%%%%%%%%%%%%%%%%%%%%%%%%%%%%%%%%%%%%%%%%%%%%%%%%%%%%%%%%%%%%%%%%
%%%%%%%%%%%%%%%%%%%%%%%%%%%%%%%%%%%%%%%%%%%%%%%%%%%%%%%%%%%%%%%%%%%%%%%%%%%%%%%%

\section{Formal multiplications and bialgebras of distributions}

In what follows all coalgebras are always assumed to be
cocommutative. We refer to \cite{Ab80} for the basics on coalgebras.

%%%%%%%%%%
\subsection{Formal maps}
%%%%%%%%%%

Let $V$ be a vector space over a field $k$ of characteristic zero.
We shall write $k[V]_i$ for the $i$th symmetric power of $V$ and
$k[V]$ for the symmetric algebra of $V$. Recall that the space
$k[V]$ is also a coalgebra: the coproduct $$\Delta: k[V]\to
k[V]\otimes k[V]$$ is
%given by
%$$\Delta(\bx)=\bx\otimes 1+ 1\otimes \bx$$ for any $\bx\in V$,
defined by the condition that all elements of $V$ are primitive, and
the counit $\epsilon: k[V]\to k$ sends an element of $k[V]$ to its
degree 0 component. We denote by $\pi_V$ the projection of $k[V]$
onto its primitive part $k[V]_1=V$.

%We will identify $k[V\times V]$ with $k[V]\otimes k[V]$.

Elements of the dual space $k[V]^*$ will be referred to as {\em
formal functions on $V$}, and those of $k[V]$ as {\em formal
distributions on $V$}. A {\em formal map $\theta$ from $V$ to $W$}
is a linear map $$\theta:k[V]\to W$$ with  $\theta(1) = 0.$

\begin{proposition}
\label{prop:morphism_prime} Any formal map $\theta \colon k[V] \to
W$ induces a unique coalgebra morphism $\theta': k[V]\to k[W]$ with
$\pi_W \theta' = \theta.$
\end{proposition}
\begin{proof}
Define the coalgebra morphism $\theta'$ by
$$\theta'(\mu) =
\sum_{n=0}^\infty \frac{1}{n!} \theta(\mu_{(1)})\cdots
\theta(\mu_{(n)}) = \epsilon(\mu) 1 + \theta(\mu) + \cdots$$ and
observe that by \cite[Corollary 2.4.17 (i)]{Ab80} any coalgebra
morphism $\theta'$ from $k[V]$ to $k[W]$ is determined by its
projection $\pi_W \theta'$.
\end{proof}

The algebra $k[V_1\times\cdots\times V_n]$ is canonically isomorphic
to $k[V_1]\otimes\cdots\otimes k[V_n]$. In order to work with formal
maps from products of vector spaces the following notation will be
of help.

The map $\pi_{V_i}:k[V_i]\to V_i$ will be denoted by $\bx_i$ and the
null map $k[V_i] \to V_i$ for any $i$ will be denoted simply by
$\bzero$ (the absence of the index $i$ should not lead to
confusion). The induced coalgebra morphism $\bx_i'$ is the identity
map on $k[V_i]$, and $\bzero' (\mu)= \epsilon(\mu) 1$. Given a
formal map
$$G\colon k[V_1\times \cdots \times V_n] \to W$$ and formal maps
$\theta_{i} \colon k[U_{i}] \to V_{i}$ for $1\leq i\leq n$ we write
$G(\theta_{1},\ldots,\theta_{n})$ for the map $G\circ \theta'_{1}
\otimes \cdots  \otimes \theta'_{n}.$

With this notation the $\bx_i$ can be treated as variables. In
particular, $G$ can be also written as $G(\bx_{1},\ldots,\bx_{n})$.
If $$G(%\bx_{1},
\ldots, \bx_{i-1},\bx_i,\bx_{i+1},
\ldots%,\bx_{n}
)=G(%\bx_{1},
\ldots, \bx_{i-1},\bzero,\bx_{i+1},
\ldots%,\bx_{n}
)$$ we say that $G$ {does not depend} on $\bx_i$ and
omit this variable altogether; the domain of definition of $G$ will
always be clear from the context.

If $V_1= \cdots =V_n=V$ the notation $G(\bx,\ldots, \bx)$ stands for
the composition of $G$ with the map $k[V]\to k[V_1\times \cdots
\times V_n]$ induced by the diagonal $V\to V_1\times\cdots\times
V_n$:
$$
\mu \mapsto \sum G(\mu_{(1)},\dots, \mu_{(n)}).
$$
Similarly one defines $G(\bx_{i_1}, \ldots, \bx_{i_n})$ when there
are various groups of repeating indices among the $i_k$.

%%%%%%%%%%
\subsection{Formal multiplications}
%%%%%%%%%%
Formal multiplications are a special form of formal maps.
\begin{definition}
A formal multiplication  on $V$ is a
formal map $$F:k[V\times V]\to V.$$
A formal multiplication on $V$ is said to be unital (or a formal loop) if
$$F|_{k[V]\otimes 1}={\pi_V}=F|_{1\otimes k[V]}.$$
\end{definition}

Since
%\begin{align*}F\in \Hom(k[V]\otimes k[V],V) &= \Hom(\bigoplus_{i,j=0}^\infty k[V]_i \otimes  k[V]_j , V)
%\\ &\cong \prod_{i,j=0}^\infty \Hom(k[V]_i\otimes k[V]_j, V)
%\end{align*}
$$F\in \Hom(k[V]\otimes k[V],V) \cong \prod_{i,j=0}^\infty \Hom(k[V]_i\otimes k[V]_j, V)$$ with our
notation we can write any unital formal multiplication $F$ as an
infinite formal sum
$$
F(\bx,\by) = \bx + \by + \sum_{i,j\geq 1} F_{i,j}(\bx,\by)
$$
with $F_{i,j}\in \Hom(k[V]_i \otimes k[V]_j, V)$, or equivalently
$$F(\mu_1 \otimes \mu_2 ) = \pi_V(\mu_1)\epsilon(\mu_2) +
\epsilon(\mu_1)\pi_V(\mu_2) + \sum_{i,j\geq 1} F_{i,j}(\mu_1 \otimes
\mu_2).$$ We say that $F_{i,j}(\bx,\by)$ is of degree $i$ in $\bx$
and $j$ in $\by$. Sometimes we shall write $\bx\by$ for a unital
formal multiplication $F(\bx,\by)$.

\begin{proposition}
\label{prop:formaldivisions}
Let $F(\bx, \by)=\bx \by$ be a unital formal multiplication. There
exist formal multiplications $\bx\backslash \by$ and $\bx/ \by$ such
that
\begin{enumerate}
\item $\bx\backslash (\bx \by) = \by = \bx( \bx \backslash \by)$ and
\item $(\by \bx)/\bx = \by = (\by/ \bx)\bx$.
\end{enumerate}
\end{proposition}

\begin{proof}
Write $F(\bx, \by) = \bx +  \by + \sum_{i,j\geq 1}
F_{i,j}(\bx, \by)$. Given a formal multiplication $D(\bx, \by)$ we have
that
\begin{displaymath}
F(\bx,D(\bx,\by))= \by \text{ if and only if } D(\bx, \by) = \by-
\bx - \sum_{i,j\geq 1} F_{i,j}(\bx,D(\bx,\by)).
\end{displaymath}
The latter recurrence determines a unique solution $D(\bx,\by)$ that
in addition satisfies $D(\bzero,\by)= F(\bzero,D(\bzero, \by)) =
\by$. By the same argument, there exists then a unique solution
$H(\bx, \by)$ to the equation $D(\bx,H(\bx,\by))= \by$. By
construction
\begin{displaymath}
H(\bx,\by) = F(\bx,D(\bx,H(\bx,\by))) = F(\bx,\by)
\end{displaymath}
so $\bx\backslash \by = D(\bx, \by)$ satisfies the required
conditions. In a similar way one proves the existence of $\bx/\by$.
\end{proof}

By Proposition~\ref{prop:morphism_prime}, any unital formal
multiplication $F(\bx,\by)$ induces a product
$$F':k[V]\otimes k[V]\to k[V].$$
Whenever we consider $k[V]$ as an algebra with multiplication
induced by $F$ we shall denote it by $k[F]$; similarly, we shall
write $k[F]^*$ for $k[V]^*$. The unit $\delta \colon k \to k[F]$ is
defined by $\alpha\mapsto \alpha 1$.

\begin{proposition}
$(k[F],\Delta,\epsilon,F',\delta)$ is an irreducible unital
bialgebra.
\end{proposition}
\begin{proof}
By construction $$F'(\mu,1) = \sum_{n=0}^\infty \frac{1}{n!}
\pi_V(\mu_{(1)})\cdots \pi_V(\mu_{(n)}) =  \mu = F'(1,\mu)$$ for any
$\mu\in k[F]$. Since $F'$ is a coalgebra morphism, the proposition
follows.
\end{proof}

In what follows we shall assume that all bialgebras in question are
unital.

%%%%%%%%%%
\subsection{The equivalence of categories}
%%%%%%%%%%
The category of unital formal multiplications is, in fact,
equivalent to that of irreducible bialgebras.

\begin{definition}
Let $F$ and $H$ be unital formal multiplications on vector spaces
$V$ and $W$ respectively. A formal map $\theta$ from $V$ to $W$ is
called a homomorphism from $F$ to $H$ (notation: $\theta\colon F
\rightarrow H$) if
\begin{displaymath}
H(\theta(\bx), \theta(\by)) = \theta(F(\bx,\by)) %\quad \bx,\by \in V
\end{displaymath}
or, equivalently, $H(\theta'(\mu_1)\otimes \theta'(\mu_2)) = \theta(
F'(\mu_1\otimes \mu_2))$ for any $\mu_1,\mu_2 \in k[V]$.
\end{definition}

Homomorphisms are the morphisms in the category of unital formal
multiplications. It follows directly from the definitions that a
homomorphism of formal multiplications $\theta\colon F\rightarrow H$
induces a homomorphism of bialgebras
$$
\theta'\colon k[F] \rightarrow k[H].
$$

\begin{proposition}
The category of unital formal multiplications and the category of
irreducible unital bialgebras are equivalent.
\end{proposition}
\begin{proof}

First, let us show that any irreducible bialgebra is isomorphic to
$k[F]$ for some unital formal multiplication $F$. The coalgebra
structure of such bialgebras is very well known: in \cite{SU02} it
is proved (see Theorem~3.2) that every such bialgebra is isomorphic
as a coalgebra to $k[V]$ for some $V$.

Given a product $m:k[V]\otimes k[V]\to k[V]$, its primitive (that
is, linear) part
$$k[V]\otimes k[V]\stackrel{m}{\longrightarrow} k[V]\stackrel{\pi_V}{\longrightarrow} V$$
is a formal multiplication $F_m$. The fact that $m(x,1)=m(1,x)=x$
means that $F_m$ is unital. It follows from the construction that
$k[F_m]$ coincides with the bialgebra $k[V]$ equipped with the
product $m$.

Now, if $\psi:k[F]\to k[H]$ is a homomorphism of bialgebras, its
primitive part $\theta=\pi_V \psi$ is a homomorphism $F\to H$ and,
clearly, $\psi=\theta'$.
\end{proof}

%%%%%%%%%%
\subsection{Analytic loops and formal loops}
\label{subsec:loopsloopsloops}
%%%%%%%%%%
Let $V$ be a finite-dimensional vector space with a basis
$\{e_1,\dots,e_n\}$ and let $\{x_1,\dots,x_n\}$ be the dual basis of
$V^*$. The symmetric algebra $k[V]$ is spanned by monomials in the
$e_i$. The dual space $k[V]^*$ is also an algebra with the {\em
convolution} product: for $f,g\in k[V]^*$ it is defined as
$$(f*g)(m)=\sum_{m_1m_2=m} f(m_1)g(m_2)$$
where $m,m_1,m_2$ are basis monomials of $k[V]$. This algebra can be
identified with the algebra $k[[x_1,\dots,x_n]]$ of formal power
series in the $x_i$, where the products of the $x_i$ are thought of
as convolutions.

Elements of $k[V]$ can be understood as elements in
$(k[V]^*)^*=k[[x_1,\dots, x_n]]^*$. For any monomial
$x_1^{a'_1}\cdots x_n^{a'_n}$ we have that $$e_1^{a_1}\cdots
e_n^{a_n}(x_1^{a'_1}\cdots x_n^{a'_n}) = x_1^{a'_1}\cdots
x_n^{a'_n}(e_1^{a_1}\cdots e_n^{a_n})$$ which, by definition of the
convolution product equals to $a_1!\cdots a_n!$ if the ordered set
of exponents $(a_1,\dots,a_n)$ coincides with $(a'_1,\dots, a'_n)$
and $0$ otherwise. Therefore $$e_1^{a_1}\cdots e_n^{a_n} =
\partial_1^{a_1}\cdots \partial^{a_n}_n\vert_{(0,\dots,0)}$$
and $k[V]$ is the coalgebra of all distributions (linear functionals
on analytic functions) on $V$ with support at zero. In particular,
the constant polynomial $1$ corresponds to the evaluation at $0$
(also known as the Dirac delta).

Now, let $G$ be an analytic local loop defined on a neighborhood of
the origin (which plays the role of the unit) in $V$. Having chosen
a basis in $V$, we may write the product $F(x,y)$ on $G$ as an
$n$-tuple
$$(F_1(x,y),\dots, F_n(x,y))$$ of power series in $2n$ variables, that
is, an element of $V\otimes k[V \times V]^*$, satisfying $F(x,0) =
x$ and $F(0,y) = y.$ Under the natural isomorphism $$V\otimes
k[V\times V]^*\cong \Hom(k[V\times V],V)$$ this condition is
equivalent to $F\vert_{k[V]\otimes 1} = \pi_V = F\vert_{1\otimes
k[V]}$. Therefore, analytic local loops give rise to formal
multiplications.

The product $F$ on the analytic local loop $G$ induces a product on
$k[V]$ which sends $\mu_1\otimes \mu_2$ to $\mu_1\cdot \mu_2 \colon
f \mapsto \mu_1\otimes \mu_2(f\circ F)$. This product is a coalgebra
map, and it gives $k[V]$ the structure of an irreducible unital
bialgebra. Since $$\mu_1\otimes \mu_2(x_i\circ F)=\mu_1\otimes \mu_2
(F_i),$$ the primitive part of $\mu_1\cdot \mu_2$ is $F(\mu_1
\otimes \mu_2)$ and, as a consequence, $$\mu_1\cdot \mu_2 = F'(\mu_1
\otimes \mu_2).$$ Our definition of the bialgebra of distributions
corresponding to a formal loop is motivated by this observation.

More generally, any analytic map $\theta \colon V \to W$ defined on
a neighbourhood of $0$ and such that $\theta(0)=0$ induces a
coalgebra morphism $\theta'$ on distributions given by $\theta'(\mu)
(f) = \mu(f \circ \theta)$ for any analytic function $f$ and any
distribution $\mu$. Note that $\theta$ gives rise to a formal map
from $V$ to $W$; if the distributions on $V$ are identified with
$k[V]$ and formal power series with $k[V]^*$, this formal map
induces the same map as $\theta'$.
%Usual
%manipulations of analytic maps such as composition or evaluation of
%some arguments in $0$ have their counterparts on formal maps and
%have been introduced.
%Finally, the map $\theta^\bullet$ corresponds to the map that the analytic map $\theta$ induces on analytic functions.

%%%%%%%%%%%%%%%%%%%%%%%%%%%%%%%%%%%%%%%%%%%%%%%%%%%%%%%%%%%%%%%%%%%%%%%%%%%%%%%%
%%%%%%%%%%%%%%%%%%%%%%%%%%%%%%%%%%%%%%%%%%%%%%%%%%%%%%%%%%%%%%%%%%%%%%%%%%%%%%%%
%%%%%%%%%%%%%%%%                                                %%%%%%%%%%%%%%%%
%%%%%%%%%%%%%%%%                                                %%%%%%%%%%%%%%%%
%%%%%%%%%%%%%%%%                   IDENTITIES                   %%%%%%%%%%%%%%%%
%%%%%%%%%%%%%%%%                                                %%%%%%%%%%%%%%%%
%%%%%%%%%%%%%%%%                                                %%%%%%%%%%%%%%%%
%%%%%%%%%%%%%%%%%%%%%%%%%%%%%%%%%%%%%%%%%%%%%%%%%%%%%%%%%%%%%%%%%%%%%%%%%%%%%%%%
%%%%%%%%%%%%%%%%%%%%%%%%%%%%%%%%%%%%%%%%%%%%%%%%%%%%%%%%%%%%%%%%%%%%%%%%%%%%%%%%

\section{Identities}

\subsection{Identities in formal loops and in bialgebras}

A {\em formal group} $F\colon k[V]\otimes k[V] \to V$ is a formal
loop that in addition satisfies the identity $F(F(\bx,\by),\bz) =
F(\bx, F(\by, \bz))$.  The consequence of this identity is that
$F'\circ (F'\otimes \Id) = F' \circ (\Id \otimes F')$ which implies
that $k[F]$ is associative. In our approach groups do not play any
special role and the bialgebras of distributions considered here
will be non-associative in general. However, the principle that
identities on loops produce identities on distributions works in
general and provides interesting examples of identities in
non-associative bialgebras.

Consider the set  of formal maps
\begin{equation}\label{eq:series}
\bx \by=F(\bx, \by), \quad  \bx\backslash \by, \quad  \bx/ \by, \quad
\bx \backslash \be,\quad \text{ and }\quad \be/ \bx
\end{equation}
with $\be=\bzero$ and $F(\bx,\by)$ is a unital formal multiplication
on a vector space $V$ (since $\bx \backslash \be$ and $\be / \bx$ do
not depend on the variable $\by$ we may consider them as defined on
$k[V]$). We may compose these formal maps to obtain new formal maps,
such as $F(\bx, \bx\backslash \be)$ (which is equal to $\be$),
$\be/(\bx\backslash \be)$ (equals to $\bx$),
$F(\bx,F(\by,F(\bx,\bz)))$, $ F(F(F(\bx,\by),\bx),\bz)$, and so on.

Let $x_1,\dots,x_n$ be a set of generators for the free loop on $n$
letters, and let $V_1,\ldots, V_n$ be $n$ copies of $V$. Since
$F(\bx,\bzero) = \bx = F(\bzero,\bx)$, by
Proposition~\ref{prop:formaldivisions} to each word
$w(x_1,\dots,x_n)$ we can assign a formal map
\begin{equation}\label{eq:w}
w:k[V_1]\otimes\ldots\otimes k[V_n]\to V
\end{equation}
by substituting $\bx_i$ for each occurrence of $x_i$ in
$w(x_1,\dots,x_n)$, and understanding the products and divisions as
in (\ref{eq:series}) above. For instance, to the word $x(y(xz))$ we
assign the map $\bx(\by(\bx\bz)) = F(\bx,F(\by,F(\bx,\bz)))$.

\begin{definition}
Given two words $u(x_1,\dots,x_n)$ and $v(x_1,\dots,x_n)$  in the
free loop on $n$ letters, we say that $F$ satisfies the identity
$$u(x_1,\dots,x_n)\sim v(x_1,\dots,x_n)$$ if the maps
$u(\bx_1,\dots,\bx_n)$ and $ v(\bx_1,\dots,\bx_n)$ coincide.
%By abuse of notation we also say that $F$ satisfies the identity $u(\bx_1,\dots,\bx_n) = v(\bx_1,\dots,\bx_n)$.
\end{definition}
Any map $w$ as in (\ref{eq:w}) induces $w'\colon
k[V_1]\otimes\cdots\otimes k[V_n] \rightarrow k[V]$. The map $w$ can
be recovered from $w'$ by taking the primitive part. In particular,
a untial formal multiplication $F$ satisfies the identity
$u(x_1,\dots,x_n)\sim v(x_1,\dots,x_n)$ if and only if
$u(\bx_1,\dots,\bx_n)'=v(\bx_1,\dots,\bx_n)'.$

Operations $\bx\backslash \by$ and $\bx/ \by$ induce the
corresponding operations on distributions; these operations were
first considered in \cite{Pe07}. We shall simply write $\mu
\backslash \nu$ and $\mu/\nu$ to denote these operations. Since  any
formal unital multiplication $\bx \by = F(\bx,\by)$ satisfies
$\bx\backslash (\bx\by)=\by =\bx(\bx\backslash \by)$ and
$(\bx\by)/\by=\bx=(\bx/\by)\by$ we have that for any $\mu,\nu\in
k[F]$
\begin{eqnarray*}
&\sum \mu_{(1)}\backslash(\mu_{(2)}\nu) = \epsilon(\mu)\nu = \sum
\mu_{(1)}(\mu_{(2)}\backslash \nu)&\\
&\sum (\mu \nu_{(1)})/\nu_{(2)}=\epsilon(\nu)\mu = \sum
(\mu/\nu_{(1)})\nu_{(2)}.&
\end{eqnarray*}
%\noindent Operations $e/x$ and $x\backslash e$ induce $\delta/\mu$ and $\mu\backslash\delta$ on
%distributions.

A particular choice of $u$ and $v$ is
$$u(x,y,z)=
x(y(x z))$$ and
$$v(x,y,z)=
((x y)x)z.$$ The identity
$u\sim v$ is called the {\em Moufang
identity}. The corresponding $u'$ and $v'$ are $u'(\mu,\nu,\eta)=\sum
\mu_{(1)}(\nu(\mu_{(2)}\eta))$ and $v'(\mu,\nu,\eta)=\sum
((\mu_{(1)}\nu)\mu_{(2)})\eta$. This shows that $F$ is a formal
Moufang loop if and only if
\begin{equation}\label{eq:linmoufang}
\sum
\mu_{(1)}(\nu(\mu_{(2)}\eta))= \sum ((\mu_{(1)}\nu)\mu_{(2)})\eta
\end{equation}
in $k[F]$.

In \cite{Pe07} $w'$ was called the {\em  linearization} of $w$. It
was proved that certain bialgebras constructed from Malcev algebras
satisfy the identity (\ref{eq:linmoufang}) above. That was
surprising since the construction of those bialgebras \cite{PS04}
had no relation with Moufang loops. Distributions provide a natural
connection between identities in loops and identities in bialgebras.
In case that we consider a local analytic loop $G$,  any word
$u(x_1,\dots,x_n)$ induces by evaluation a map $u\colon G\times
\cdots \times G \to G$ and  a coalgebra map $u'$ on distributions
which agrees with the map $u'$ defined above. Therefore, $G$
satisfies the identity $u\sim v$ in the usual sense if and only if
the formal loop corresponding to $G$ satisfies the identity $u\sim
v$.

%%%%%%%%%%%%%%%%%%%%%%%%%%%%%%%%%%%%%%%%%%%%%%%
\subsection{Right alternativity}\label{subsec:mono}
%%%%%%%%%%%%%%%%%%%%%%%%%%%%%%%%%%%%%%%%%%%%
Another example of an identity on formal multiplications is
$$x(yy)=(xy)y$$
called {\em right alternativity}. The corresponding bialgebra
identity reads
\begin{equation}\label{eq:linralt}
\sum \mu(\nu_{(1)}\nu_{(2)})=
\sum (\mu\nu_{(1)})\nu_{(2)}.
\end{equation}

It was proved by Sabinin and Mikheev \cite{SM85} that the right
alternativity in a formal loop implies the identity
\begin{equation}\label{eq:monoalt}
x(y^ky^l)=(xy^{k})y^{l}
\end{equation}
for all $k,l\geq 0$ (see also \cite{S99}).

The importance of the right alternativity for Lie theory was
understood first by Sabinin and Mikheev \cite{SM87, MS90}. They
realized that this algebraic property for a local loop is satisfied
if and only if the loop comes from a flat connection as a so-called
{\em geodesic loop}, and showed that, in fact, the multiplication in
{any} local loop can be modified so as to become right alternative.

Given a local loop $(G, \cdot)$ with the unit $e$, one can define
the {\em canonical} flat connection $\nabla$ on the tangent bundle
to $G$ in a neighbourhood of $e$ as follows. For $a,b\in G$ two
points in a small neighbourhood of $e$, the parallel transport of
the tangent space $\mathrm{T}_bG$ to $\mathrm{T}_aG$ is induced by a
self-map of $G$ that sends $x$ to $a\cdot(b\backslash x)$. There is
a new product on $G$ given by
\begin{equation}\label{eq:monoaltproduct}
a\times b =\exp_a{(a\log{b})},
\end{equation}
where $\exp_a$ is the exponential map of the connection $\nabla$ at
the point $a$ and $\log_a$ is its inverse (we write simply $\exp$
and $\log$ for $\exp_e$ and $\log_e$ respectively) and $a\log{b}$
stands for the parallel transport of the vector $\log{b}\in
\mathrm{T}_eG$ to $\mathrm{T}_aG$. The local loop $(G,\times)$ is
then right alternative. Note that the canonical connection for the
loop $(G,\times)$  is the same as that of $(G,\cdot)$. The original
local loop $(G,\cdot)$ can be reconstructed from $(G,\times)$ and
the operation $\Phi(a,b)$ defined by
\begin{equation}\label{eq:phi}
a\times\Phi(a,b)=a\cdot b.
\end{equation}

The right alternative modification can also be defined for formal
loops. Given a formal loop $F$ on a vector space $V$, the {\em
formal canonical connection} of $F$ is the restriction of $F$ to the
subspace $$k[V]\otimes V\subset k[V]\otimes k[V].$$ Let us say that
two formal loops on the same vector space are {\em similar} if their
formal canonical connections coincide.
\begin{lemma}\label{lem:similarloop}
Each formal loop is similar to a unique right alternative formal loop.
\end{lemma}
\begin{proof}
A unital formal multiplication can be written as
$$F(\bx,\by)=\bx + \by + q_{1}(\bx,\by)+q_{2}(\bx, \by)+\ldots$$
where $q_{j}(\bx, \by) = \sum_{i=1}^\infty F_{i,j}(\bx,\by)$.
Specifying the canonical connection for $F$ is
the same thing as specifying $q_1(\bx,\by)$.

Given  $q_1(\bx,\by)$, the  right alternative formal loop similar to
$F$ can be reconstructed inductively. Assume that the $q_i(\bx,\by)$
with $i<n$ are known, and consider the equation
$F(F(\bx,\by),\by)=F(\bx,F(\by,\by))$, modulo the terms of degree
$>n$ in $\by$. A simple calculation shows that, apart from the
(compositions of the) $q_i$ with $i<n$, this equation contains the
term $q_n(\bx,\by)$ with coefficient $2$ on the left-hand side and
$2^n$ on the right-hand side. Therefore, for $n>1$ we see that $q_n$
can be expressed via the $q_i$ with $i<n$.
\end{proof}

Let $\Phi:k[V]\otimes k[V]\to V$ be a formal multiplication such
that $\Phi|_{1\otimes k[V]}=\pi_V$ and $\Phi|_{k[V]_{\geq 1}\otimes
(1\oplus V)}=0$. In other words,
$$\Phi(\bx,\by)= \by+\sum_{i\geq 1, j\geq 2}\Phi_{i,j}(\bx, \by),$$
with $\Phi_{i,j}(\bx, \by)$ of degree $i$ in $\bx$ and of degree $j$ in $\by$. Call such a
multiplication a {\em similarity}.

\begin{lemma}\label{lem:similarity}
Two formal loops $F_1$ and $F_2$ are similar if and only if there is
a similarity $\Phi$ such that $F_1(\bx,\by)=F_2(\bx,\Phi(\bx,\by))$.
\end{lemma}
The proof is straightforward. For the ``if'' part compare the
corresponding homogeneous terms $(F_1)_{i,j}$ and $(F_2)_{i,j}$; for
the ``only if'' part define $\Phi$ inductively by the degree. \qed

\medskip

These notions have their versions for bialgebras. If $\Phi$ is a
similarity, we obtain a coalgebra morphism
$\Phi':k[V]\otimes k[V]\to k[V]$, which satisfies
\begin{equation}
\label{eq:similin} \Phi'(\mu_1,1) =\epsilon(\mu_1)1, \enskip
\Phi'(1,\mu_2) = \mu_2 \enskip \text{and} \enskip \Phi'(\mu_1,\alpha) =\epsilon(\mu_1)\alpha
\end{equation}
for all $\mu_1, \mu_2 \in k[V]$ and $\alpha \in V$. Conversely, the primitive part
of a coalgebra morphism $\Phi'$ satisfying these conditions is a
similarity.

Let $F_1$ and $F_2$ be two similar formal loops with
$F_1(\bx,\Phi(\bx,\by))=F_2(\bx,\by)$. Denote by $\times$ and $\cdot$ and the
products in $k[F_1]$ and in $k[F_2]$ respectively. Then we have
\begin{equation}
\label{eq:linphi}
\sum \mu_{(1)}\times \Phi'(\mu_{(2)},\nu)=\mu\cdot\nu.
\end{equation}
If $k[V]$ has two different bialgebra products $\times$ and $\cdot$
such that there exists a map $\Phi'$ satisfying (\ref{eq:similin})
and (\ref{eq:linphi}), we say that the products $\times$ and $\cdot$
are {\em similar} and that $\Phi'$ is a {\em (bialgebra) similarity}
between them.
\begin{lemma}\label{lem:simprod}
If $\times$ and $\cdot$ are two similar products on $k[V]$, then for
any $\mu\in k[V]$ and $\alpha\in V$
$$\mu\times\alpha=\mu\cdot\alpha.$$
\end{lemma}
The proof is an immediate consequence of (\ref{eq:similin}) and
(\ref{eq:linphi}).

%%%%%%%%%%%%%%%%%%%%%%%%%%%%%%%%%%%%%%%%%%%%%%%%%%%%%%%%%%%%%%%%%%%%%%%%%%%%%%%%
%%%%%%%%%%%%%%%%%%%%%%%%%%%%%%%%%%%%%%%%%%%%%%%%%%%%%%%%%%%%%%%%%%%%%%%%%%%%%%%%
%%%%%%%%%%%%%%%%                                                %%%%%%%%%%%%%%%%
%%%%%%%%%%%%%%%%                                                %%%%%%%%%%%%%%%%
%%%%%%%%%%%%%%%%         BIALGEBRAS OF DISTIBUTIONS AND         %%%%%%%%%%%%%%%%
%%%%%%%%%%%%%%%%                SABININ ALGEBRAS                %%%%%%%%%%%%%%%%
%%%%%%%%%%%%%%%%                                                %%%%%%%%%%%%%%%%
%%%%%%%%%%%%%%%%%%%%%%%%%%%%%%%%%%%%%%%%%%%%%%%%%%%%%%%%%%%%%%%%%%%%%%%%%%%%%%%%
%%%%%%%%%%%%%%%%%%%%%%%%%%%%%%%%%%%%%%%%%%%%%%%%%%%%%%%%%%%%%%%%%%%%%%%%%%%%%%%%

%%%%%%%%%
\section{Bialgebras of distributions and Sabinin algebras}

%%%%%%%%%
\subsection{Shestakov--Umirbaev's functor $\SU$}\label{subsec:shu}
%%%%%%%%%

Let $\mathcal{S}$ be a set. Denote by $k\{\mathcal{S}\}$  the unital
free non-associative algebra generated by the elements of
$\mathcal{S}$. The algebra $k\{\mathcal{S}\}$ can be given a
structure of an irreducible bialgebra: the comultiplication is
defined by the condition that all the elements of $\mathcal{S}$ are
primitive; the counit $\epsilon \colon k\{\mathcal{S}\}\rightarrow
k$ is the homomorphism that sends $1$ to $1$ and all the elements of
$\mathcal{S}$ to $0$.

Recall that instead of the antipodes, non-associative bialgebras
have operations of left and right division $\backslash$ and $/$. In
$k\{\mathcal{S}\}$ they are as follows. Starting with $1\backslash v
= v$, $a\backslash v = -a v$ for any generator $a\in \mathcal{S}$
and $v\in k\{\mathcal{S}\}$ we use induction on the degree $\vert
u\vert$ of of $u$ to define a bilinear operation $u\backslash v$ so
that
\[
\sum u_{(1)} \backslash (u_{(2)}v) =\epsilon(u) v.
\]
We also have $\sum u_{(1)}(u_{(2)}\backslash v)
=\epsilon(u)v$. Indeed, by induction on $\vert u\vert$, we get
\begin{multline*}
uv=\sum {u_{(1)}}\left({u_{(2)}}\backslash
u_{(3)}v\right)\hskip -1pt =\sum u_{(1)}(u_{(2)}\backslash v)
+\sum_{\vert u_{(1)}\vert< \vert u\vert} \epsilon(u_{(1)}) u_{(2)} v
\\=-\epsilon(u)v+uv +\sum u_{(1)}(u_{(2)}\backslash v).
\end{multline*}
Hence $\sum
u_{(1)}(u_{(2)}\backslash v) =\epsilon(u)v$. Similarly we define a
bilinear operation $u/v$ that satisfies $$\sum (u v_{(1)})v_{(2)}
=\epsilon(v)u = \sum (u/v_{(1)})v_{(2)}.$$

Apart from the generators and their linear combinations, the algebra
$k\{\mathcal{S}\}$ has many other primitive elements. All these
elements were described by Shestakov and Umirbaev in \cite{SU02}.

Let $u=((x_1x_2)\cdots)x_m$ and $v=((y_1y_2)\cdots )y_n$ with $x_i$
and $y_j$ primitive. The {\em primitive operations} $p(u,v,z)$ are
defined by
\[
p(x_1,\dots,x_m;y_1,\dots,y_n;z)=p(u,v,z) = \sum
(u_{(1)}v_{(1)})\backslash (u_{(2)},v_{(2)},z)
\]
where $(x,y,z) = (xy)z-x(yz)$ denotes the associator and $z$ is
primitive. In \cite{SU02} it is shown that the $p(u,v,z)$ are
primitive, and, moreover, that each primitive element of
$k\{\mathcal{S}\}$ can be obtained from the generators by applying
repeatedly the commutators and the operations $p(u,v,z)$, and taking
linear combinations.

Since $p(x_1,\dots,x_m;y_1,\dots,y_n;z)$ are just polynomial
expressions in $x_1,\dots,x_m$, $y_1,\dots,y_n$, $z$, they make
sense in any algebra; given a non-associative algebra $A$ we shall
consider them as new $(m+n+1)$--ary operations obtained from the
product on $A$. When evaluating in an arbitrary algebra $A$, the
compact notation $p(u,v,z)$ for the operation
$p(x_1,\dots,x_m;y_1,\dots,y_n;z)$ may be misleading since it
suggests that we should first evaluate $u=((x_1x_2)\cdots)x_m$ and
$v=((y_1y_2)\cdots )y_n$ and then apply a ternary operation
$p(u,v,z)$. In order to avoid confusion, we shall write
$p(\underline{u},\underline{v},z)$ when working in a
non-associative, not necessarily free, algebra $A$. The relation
\[
(u,v,z) = \sum u_{(1)}v_{(1)}\cdot
p(\underline{u}_{(2)},\underline{v}_{(2)},z)
\]
also makes sense in any algebra $A$ even if it is not a bialgebra.
This is a consequence of the corresponding identity in
$k\{\mathcal{S}\}$ . Observe that when $u, v, u_{(1)}$ and $v_{(1)}$
are not used as the arguments of $p$, they become products in $A$,
so we do not need to underline them.

Shestakov and Umirbaev related their work with the results of
Mikheev and Sabinin on local loops \cite{SM87, MS90}. Namely, in
\cite{SU02} they defined for any non-associative algebra $A$ the
operations
\begin{eqnarray*}
\langle 1;y,z\rangle &=& \langle y,z\rangle = -[y,z] = -yz + zy\\
\langle x_1,\dots,x_m;y,z\rangle &=& \langle
\underline{u};y,z\rangle =  -p(\underline{u},y,z)
+p(\underline{u},z,y)\\
\Phi^{SU}(x_1,\dots,x_m;y_1,\dots,y_n) &=&\\
&&\hskip -3cm  \frac{1}{m!}\frac{1}{n!}\sum_{\tau\in S_m,\sigma\in
S_n} p(x_{\tau(1)},\dots,x_{\tau(m)};y_{\sigma(1)},\dots,
y_{\sigma(n-1)};y_{\sigma(n)})
\end{eqnarray*}
with $u=((x_1x_2)\cdots)x_m$, $S_m$ the symmetric group on $m$
letters and $m\geq 1, n\geq 2$. With these operations $A$ turns out
to be a Sabinin algebra \cite{SU02} so we have a functor from
non-associative algebras to Sabinin algebras
\[
A\mapsto \SU(A)
\]
that generalizes the usual functor from associative algebras to Lie
algebras given by assigning to an associative algebra its commutator
algebra. The primitive elements of any bialgebra $W$ form a Sabinin
subalgebra of $\SU(W)$.

One is then naturally led to ask whether every Sabinin algebra is
isomorphic to a Sabinin algebra of the primitive elements in some
irreducible bialgebra. An affirmative answer (with a modified
version of the operations $p(\,;\,;\,)$ and, hence, of the functor
$\SU$) was given in \cite{Pe07}. Given a Sabinin algebra $(V,\langle
\,;\,,\,\rangle,\Phi')$ the corresponding bialgebra is denoted by
$U(V)$ and has the following universal property: any homomorphism of
Sabinin algebras from $V$ to a unital algebra $A$ extends to a
unique homomorphism of unital algebras $U(V)\rightarrow A$. The
algebra $U(V)$ was called in \cite{Pe07} {\em the universal
enveloping algebra} of $V$.

There is a Poincar\'e--Birkhoff--Witt Theorem for the universal
enveloping algebras of Sabinin algebras: as a coalgebra, $U(V)$ is
isomorphic to $k[V]$.
%{(According to Theorem 3.2 in \cite{SU02} and its proof?)}
Moreover, the algebra $U(V)$ is filtered and the
corresponding associated graded algebra is commutative and
associative: it is isomorphic to the symmetric algebra $S(V)$. If we
start with an irreducible bialgebra $W$, $\Prim(W)$ is a Sabinin
subalgebra of $\SU(W)$ and if $\{e_1,e_2,\dots,e_\alpha,\dots\}$ is
a basis of $\Prim(W)$ then
$$\{((e_{i_1}e_{i_2})\cdots )e_{i_k}\mid 0\leq i_1\leq\cdots\leq i_k,
k\geq 0\}$$ is a basis of $W$ (\emph{Poincar\'e--Birkhoff--Witt
basis}). The universal property of the enveloping algebras
gives an isomorphism
\[
U(\Prim(W))\cong W
\] of bialgebras, which identifies the respective Poincar\'e--Birkhoff--Witt
bases. In this way irreducible bialgebras can be classified in terms
of the Sabinin algebra of their primitive elements. In the sequel we shall often
write irreducible bialgebras as pairs $(k[V], \cdot)$ where $\cdot$ is a product
on the coalgebra $k[V]$. Sometimes, for clarity, we shall also indicate the product
explicitly while working with primitive operations and the bialgebra divisions.

One useful consequence of the Poincar\'e--Birkhoff--Witt Theorem is the following
\begin{lemma}\label{lem:powers}
Any irreducible bialgebra is additively spanned by elements of the
form $x^n=((xx)x\ldots )x$ with $n\geq 0$ and $x$ primitive.
\end{lemma}

%%%%%%%%%%
\subsection{Similarity of bialgebras and the primitive operations}
%%%%%%%%%%

In a Sabinin algebra the identities for the brackets do not involve
the multioperator, and vice versa. Here we shall see how to modify a
product in a bialgebra so that the bracket operations defined via
the Shestakov-Umirbaev operations do not change and so that
$\Phi^{SU}$ takes any prescribed form.

\begin{proposition}
\label{prop:invariance} Let $k[V]$ be a bialgebra with respect to
two similar products $\cdot$ and $\times$. Then for any $\mu\in
k[V]$ and $\alpha,\beta\in V$
\begin{displaymath}
\langle \underline{\mu};\alpha,\beta\rangle ^\cdot = \langle
\underline{\mu};\alpha,\beta\rangle^\times.
\end{displaymath}
\end{proposition}
\begin{proof}

Let $\Phi'$ be the similarity between $\times$ and $\cdot$. By the
definition of $p^\cdot(\underline{\mu},\underline{\nu},\alpha)$,
(\ref{eq:linphi}) and Lemma~\ref{lem:simprod}, these operations can
be written in terms of $\times$ and $\Phi'$ as
\begin{align*}
\label{eq:Pop}
p^\cdot(\underline{\mu},\underline{\nu},\alpha)&=\sum (\mu_{(1)}\cdot\nu_{(1)})\backslash^\times
(\mu_{(2)},\nu_{(2)},\alpha)^\cdot\\
&=\epsilon(\mu)\epsilon(\nu)\alpha -\sum
(\mu_{(1)}\times\Phi'(\mu_{(2)},\nu_{(1)}))\backslash^\times(\mu_{(3)}\times
\Phi'(\mu_{(4)},\nu_{(2)}\times\alpha))
\end{align*}
Hence, by (\ref{eq:similin}) we have
\begin{eqnarray*}
 p^\cdot(\underline{\mu},\alpha,\beta) &=& -\sum (\mu_{(1)}\times
\Phi'(\mu_{(2)},\alpha_{(1)}))\backslash^\times (\mu_{(3)}\times
\Phi'(\mu_{(4)},\alpha_{(2)}\times \beta))\\ & =&-\sum
(\mu_{(1)}\times \alpha )\backslash^\times (\mu_{(2)}\times \beta)
\\
&&-\sum \mu_{(1)}\backslash^\times (\mu_{(2)}\times
\Phi'(\mu_{(3)},\alpha\times \beta))\\
&=& -\sum (\mu_{(1)}\times
\alpha)\backslash^\times (\mu_{(2)}\times\beta)
-\Phi'(\mu,\alpha\times\beta).
\end{eqnarray*}
It follows that
\begin{eqnarray*}
-p^\cdot(\underline{\mu},\alpha,\beta) +
p^\cdot(\underline{\mu},\beta,\alpha) &=& \sum (\mu_{(1)}\times
\alpha) \backslash^\times (\mu_{(2)}\times \beta)- \\ && -\sum
(\mu_{(1)}\times \beta) \backslash^\times (\mu_{(2)}\times
\alpha)+\Phi'(\mu,[\alpha,\beta]) \\ &=& \sum (\mu_{(1)}\times
\alpha) \backslash^\times (\mu_{(2)}\times \beta)- \\ &&-\sum
(\mu_{(1)}\times \beta) \backslash^\times (\mu_{(2)}\times \alpha) +
\epsilon(\mu)[\alpha,\beta],
\end{eqnarray*}
an expression that does not depend on the particular $\Phi'$.
%With $\Phi'(\mu,\nu) = \epsilon(\mu)\nu$ we get the desired equality.
\end{proof}

\begin{proposition}
\label{prop:changingthephi}
Let $k[V]$ be a bialgebra with respect to the product $\cdot$.
Given any set of multilinear operations
\[\Phi=\{\Phi_{i,j}:k[V]_i\otimes k[V]_j\to V\}\]
for $i\geq 1$ and $j\geq 2$ there exists a product $\times$ on
$k[V]$ similar to $\cdot$, such that the operations $\Phi^{SU}$ in
$(k[V], \times)$ coincide with the operations $\Phi$.
\end{proposition}

\begin{proof}
%The maps $\Phi_{i,j}$ act on the symmetric tensors of degree $i$ and $j$. - ya no hablamos de tensores.
Extend the definition of the $\Phi_{i,j}$ to the cases $i=0$ and
$j=1$ by setting $\Phi_{0,j}$ and  $\Phi_{i,1}$ to be identically
zero. Take $$\Psi(x,1)=\epsilon(x)1$$ for any $x\in k[V]$ and define
the bialgebra similarity $$\Psi\colon k[V]\otimes k[V]\rightarrow
k[V]$$ inductively by
\[
\Psi(x,b^{m+1}) = \sum x_{(1)}\backslash^\cdot
\Bigl((x_{(2)}\cdot\Psi(x_{(3)},b^m_{(1)}))\cdot\bigl(
\epsilon(x_{(4)})\epsilon(b^m_{(2)})
b-\Phi(\underline{x}_{(4)};\underline{b}^m_{(2)}, b)\bigr)\Bigr)
\]
for any $x\in k[V]$ and $b\in V$. Here $b^{m+1}$ stands for
$((b\cdot b)\ldots)\cdot b$. According to Lemma~\ref{lem:powers}
this determines $\Psi$ completely. It is easy to check that
$\Psi(x,b) = \epsilon(x)b$, and an induction on $m$ shows that
$\Psi(1,b^m) = b^m$, $\Delta(\Psi(x,b^m))=\Psi\otimes
\Psi(\Delta(x,b^m))$ and
$\epsilon(\Psi(x,b^m))=\epsilon(x)\epsilon(b^m)$.

Define a new product $\times$ on $k[V]$ by setting
\[
x\times y=\sum x_{(1)}\cdot\Psi(x_{(2)},y).
\]
In $(k[V],\times)$, on one hand,
$$(x,b^m,b)=(x_{(1)}\times b^m_{(1)})\times\Phi^{SU}(\underline{x}_{(2)};\underline{b}^m_{(2)},b),$$
and, on the other hand,
\begin{align*}
(x,b^m,b)&=\sum (x_{(1)}\cdot\Psi(x_{(2)},b^m))\cdot b -\sum x_{(1)}\cdot\Psi(x_{(2)},b^{m+1})\\
&=\sum(x_{(1)}\cdot\Psi(x_{(2)},b^m_{(1)}))\cdot\Phi(\underline{x}_{(3)};\underline{b}^m_{(2)},  b)\\
&=\sum(x_{(1)}\times b^m_{(1)})\times\Phi(\underline{x}_{(2)};\underline{b}^m_{(2)}, b).
\end{align*}
Using these two ways of computing $\sum (x_{(1)}\times b^m_{(1)})\backslash^\times
(x_{(2)},b^m_{(2)},b)$ we get $\Phi^{SU}=\Phi$ as desired.
\end{proof}

%%%%%%%%%%
\subsection{The equivalence of categories}
%%%%%%%%%%

It is known from \cite{Pe07} that the category of irreducible
bialgebras is equivalent to that of Sabinin algebras. The proof
given in \cite{Pe07}, however, uses primitive operations different
from the original operations $p(x_1,\dots,x_m;y_1,\dots,y_n;z)$
considered by Shestakov and Umirbaev. Here we shall show that the
functor that assigns to an irreducible bialgebra its subspace of
primitive elements with the operations defined in the
Section~\ref{subsec:shu} also gives an equivalence of categories.

\begin{lemma}
\label{lem:presentation} Let $W$ be an irreducible bialgebra,
$(\Prim(W),\langle\,;\,,\,\rangle',\Phi')$ the Sabinin subalgebra of
its primitive elements and
$\mathcal{V}=\{e_1,\ldots,e_{\alpha},\ldots\}$ - a basis of the
vector space $\Prim(W)$. Let $k\{\mathcal{V}\}$ be the unital free
non-associative algebra on $\mathcal{V}$ and $I$ the ideal of $k\{
\mathcal{V}\}$ generated by
\[
\langle u;a,b\rangle -\langle \underline{u};a,b\rangle'\quad \text{
and }\quad  \Phi^{SU}(u;v) -\Phi'(\underline{u},\underline{v})
\]
for any $a,b\in \Prim(W)$ and $u,v$ right-normed\footnote{that is,
of the form $((e_{i_1}e_{i_2})\dots)$, or, in other words, with all
opening brackets to the left of the first argument.} monomials in
the $e_i$. Then $W\cong k\{ \mathcal{V}\}/I$.
\end{lemma}
\begin{proof}
Denote by $\bar{W}$ be the algebra $k\{ \mathcal{V}\}/I$ and by
$\pi\colon k\{\mathcal{V}\}\rightarrow \bar{W}$ the quotient map.
Since $\mathcal{V}$ is a basis of $\Prim(W)$, there is an
epimorphism $k\{\mathcal{V}\}\rightarrow W$ which vanishes on $I$,
and, hence, factors through an epimorphism $\varphi\colon
\bar{W}\rightarrow W$. In order to show that $\varphi$ is an
isomorphism, we exhibit a vector space basis of $\bar{W}$ which is
sent by $\varphi$ to the Poincar\'e--Birkhoff--Witt basis of $W$.

By definition, in $\bar{W}$ we have
\begin{equation}
ua\cdot b -ub\cdot a = -\sum u_{(1)}\langle
\underline{u}_{(2)};a,b\rangle' \label{eq:flip}
\end{equation}
for any right-normed monomial $u$ in $\pi(\mathcal{V})$ and $a,b\in
\pi(\mathcal{V})$. It follows that any two right-normed monomials in
$\pi(\mathcal{V})$ that differ only by a permutation of their
variables, are equal in $\bar{W}$ modulo monomials of smaller
degree. Using this fact, together with the definition of $\Phi^{SU}$
we see that, modulo the right-normed monomials of lower order
\begin{equation}
uv\cdot a - u\cdot va\equiv \sum (u_{(1)} v_{(2)})
\Phi'(\underline{u}_{(2)};\underline{v}_{(2)})\equiv 0 \label{eq:product_of_powers}
\end{equation}
for any pair of right-normed monomials $u$ and $v$ in
$\pi(\mathcal{V})$ and $a\in \pi(\mathcal{V})$.

Using the induction on the degree of the monomials we now can deduce
that $\bar{W}$ admits a Poincar\'e--Birkhoff--Witt type set of
linear generators $((\bar{e}_{i_1}\bar{e}_{i_2})\cdots
)\bar{e}_{i_k}$ where $0\leq i_1\leq\cdots\leq i_k, k\geq 0$ and
$\bar{e}_{\alpha}=\pi({e}_{\alpha})$. Since $\varphi$ sends this set
to a Poincar\'e--Birkhoff--Witt basis of $W$ then it must be a basis
of $\bar{W}$ and $\bar{W}\cong W$.
\end{proof}

\begin{theorem}
The functor from the category of irreducible bialgebras to that of
Sabinin algebras, which assigns to a bialgebra $W$ the Sabinin
subalgebra $\Prim(W)$ of $\SU(W)$ is an equivalence of categories.
\end{theorem}
\begin{proof}
We will show that the functor $W\mapsto \Prim(W)$ is (1) faithful, (2)
full and (3) essentially surjective.

(1) Recall that any irreducible bialgebra, as a coalgebra, is
isomorphic to $k[V]$ where $V$ is the space of the primitive
elements. Therefore, by Corollary~2.4.17 in \cite{Ab80} any
homomorphism $W\rightarrow W'$ of bialgebras, with $W$ irreducible,
is determined by its restriction to $\Prim(W)$. This implies that
the functor $W\rightarrow \Prim(W)$ is faithful.

(2) This is a consequence of Lemma \ref{lem:presentation}.

(3) It was shown in \cite{Pe07} that given a Sabinin algebra
$(V,\langle\,;\,,\,\rangle,\Phi')$ there exists an irreducible
cocommutative unital bialgebra $(k[V],\cdot)$ such that
 the operations $\langle \,;\,,\,\rangle$ are recovered as
$$
(x \cdot a)\cdot b-(x\cdot b)\cdot a =-\sum x_{(1)}
\cdot \langle\underline{x}_{(2)};a,b\rangle.
$$
Now, by Proposition~\ref{prop:changingthephi} the product $\cdot$
can always be modified in such a way that the operations
$\langle\,;\,,\,\rangle$ remain the same and that the multioperator
on $V$ takes any desired form.
\end{proof}

%%%%%%%%%%%%%%%%%%%%%%%%%%%%%%%%%%%%%%%%%%%%%%%%%%%%%%%%%%%%%%%%%%%%%%%%%%%%%%%%
%%%%%%%%%%%%%%%%%%%%%%%%%%%%%%%%%%%%%%%%%%%%%%%%%%%%%%%%%%%%%%%%%%%%%%%%%%%%%%%%
%%%%%%%%%%%%%%%%                                                %%%%%%%%%%%%%%%%
%%%%%%%%%%%%%%%%                                                %%%%%%%%%%%%%%%%
%%%%%%%%%%%%%%%%   SABININ ALGEBRAS AND FORMAL MULTIPLICATIONS  %%%%%%%%%%%%%%%%
%%%%%%%%%%%%%%%%                                                %%%%%%%%%%%%%%%%
%%%%%%%%%%%%%%%%                                                %%%%%%%%%%%%%%%%
%%%%%%%%%%%%%%%%%%%%%%%%%%%%%%%%%%%%%%%%%%%%%%%%%%%%%%%%%%%%%%%%%%%%%%%%%%%%%%%%
%%%%%%%%%%%%%%%%%%%%%%%%%%%%%%%%%%%%%%%%%%%%%%%%%%%%%%%%%%%%%%%%%%%%%%%%%%%%%%%%

%%%%%%%%%%
\section{Sabinin algebras and formal multiplications}
%%%%%%%%%%

In this section we show directly, following the method of Sabinin
and Mikheev, that the category of Sabinin algebras and that of
unital formal multiplications are equivalent. As a result, we shall
have two constructions of a Sabinin algebra associated with a formal
multiplication: via the primitive elements in the bialgebra of
distributions, described in the preceding two sections, and the
direct construction of the present section. These two constructions,
however, do not coincide. We shall prove that the operations
$\langle\,;\,,\,\rangle$ are the same in both cases and exhibit a
formal multiplication for which the two multioperators are
different.

\subsection{The geometry of the operations in a Sabinin algebra}

For a Lie group $G$ the left multiplication by elements of $G$ gives
a flat connection (the {\em canonical connection}) on the tangent
bundle of $G$. All covariant derivatives of the torsion tensor of
the canonical connection vanish and the torsion tensor itself
coincides on the tangent space to the unit, up to sign, with the
bracket of the Lie algebra of $G$.

A generalization of this approach led Sabinin and Mikheev to the
first successful general treatment of the non-associative Lie
theory. They observed that an infinitesimal loop satisfying the
right alternative identity is, essentially, the same thing as a germ
of a flat affine connection. It is known that an (analytic) flat
affine connection can be reconstructed locally from its torsion
tensor and its covariant derivatives; therefore, these tensors
provide analogues of the Lie brackets for right alternative
infinitesimal loops. The identities for these operations are the
universal identities satisfied by the covariant derivatives of the
torsion tensor of a flat affine connection; their explicit form is
well-known.

Any infinitesimal loop determines a unique right alternative
infinitesimal loop; and the difference between the two is measured
by the operation $\Phi(a,b)$ defined by the equation (\ref{eq:phi}).
If this operation is analytic, it is reconstructed from its Taylor
series in the normal coordinates (local coordinates on the loop
coming from the tangent space via the exponential map). The
homogeneous terms of this Taylor series form a set of multilinear
operations (multioperator) which complements the derivatives of the
torsion tensor as a part of the structure of a Sabinin algebra.

These constructions can be translated into the formal setting with
minimal effort, as we shall now see.

\subsection{The torsion of a formal flat connection and the Mikheev-Sabinin brackets}

A {\em formal vector field} is a linear map $A: k[V]\to V$. The
product of a formal vector field $A$ with a formal function $f$ is
given by
$$
fA\colon \mu \mapsto \sum f\bigl(\mu_{(1)}\bigr)
A\bigl(\mu_{(2)}\bigr).
$$
This action provides the formal vector fields with the structure of
a free $k[V]^*$--module. In fact, any set $\{A_i\}_i$ of formal
vector fields such that $\{A_i(1)\}$ is a basis of $V$ gives a
$k[V]^*$--basis of $\Hom(k[V],V)$.

A formal vector field $A$ gives a derivation $D_A$ of the algebra
$k[V]^*$ of formal functions into itself:
$$
D_A(f) = A(f)\colon \mu \mapsto \sum f\bigl(\mu_{(1)}
A(\mu_{(2)})\bigr)
$$
where the product on $k[V]$ is that of the symmetric algebra. We
have $(fA)(g)= f\cdot A(g)$. Formal vector fields form a Lie algebra
with the Lie bracket $[A,B]$ given by
$$
[A,B]\colon \mu \mapsto \sum B\bigl(\mu_{(1)}A(\mu_{(2)})\bigr) -
A\bigl(\mu_{(1)} B(\mu_{(2)})\bigr).
$$
Clearly $[D_A,D_B] = D_{[A,B]}$.  We also have that $$[A,fB] = A(f)B
+ f[A,B].$$

A {\em formal flat affine connection} is a linear map $k[V]\otimes
V\to V$ whose restriction to $1\otimes V$ is the identity. For a
given formal connection, $\mu \in k[V]$ and $v\in V$, we write $\mu
* v$ for the image of $\mu \otimes v$. The
vector field $v^*\colon \mu \mapsto \mu
* v$ is said to be \emph{adapted} to the \emph{tangent vector} $v$. There always exists a
unique ``inverse'' map $k[V]\otimes V\to V$ sending $\mu \otimes u$
to an element that we denote by $\mu \backslash^* u$ and such that
$\sum \mu_{(1)}\backslash^* (\mu_{(2)}*v)= \epsilon(\mu) v =
\sum \mu_{(1)} * (\mu_{(2)} \backslash^* v)$.

The covariant differentiation with respect to the formal vector
field $A$ is defined as
$$
\nabla_A(B)\colon \mu \mapsto \sum
B\bigl(\mu_{(1)}A(\mu_{(2)})\bigr) -
\bigl(\mu_{(1)}A(\mu_{(2)})\bigr)*\bigl(\mu_{(3)}\backslash^*
B(\mu_{(4)})\bigr).
$$
\begin{proposition}
Let $A,B$ be formal vector fields, $f$ a formal function and $v,w
\in V$. Then
\begin{enumerate}
\item $\nabla_{fA}(B) = f\nabla_A(B)$,
\item $\nabla_{A}(fB) = A(f)B + f \nabla_A(B)$,
\item $\nabla_{v^*} (w^*) = 0$.
\end{enumerate}
\end{proposition}
\begin{proof}
We shall only prove (3). By definition
\begin{align*}
\nabla_{v^*}(w^*) (\mu) &= \sum w^*(\mu_{(1)}v^*(\mu_{(2)}))
- (\mu_{(1)}v^*(\mu_{(2)}))*(\mu_{(3)}\backslash^* w^*(\mu_{(4)}))
\\
&= \sum  (\mu_{(1)}(\mu_{(2)}*v))* w -
(\mu_{(1)}(\mu_{(2)}*v))*(\mu_{(3)}\backslash^* (\mu_{(4)}*w))\\
&=  \sum  (\mu_{(1)}(\mu_{(2)}*v))* w -
(\mu_{(1)}(\mu_{(2)}*v))*(\epsilon(\mu_{(3)})w) = 0.
\end{align*}
\end{proof}

The torsion of two formal vector fields $A$ and $B$ is defined in the usual way
$$
T(A,B) = \nabla_A(B) - \nabla_B(A) - [A,B].
$$
In the case of adapted vector fields $x^*, y^*$ with $x,y\in V$ we get
$$
T(x^*,y^*) = -[x^*,y^*].
$$

%For a
%given formal connection, $s\in k[V]$ and $v\in V$, we write $sv$ for
%the image of $s\otimes v$. There always exists a unique ``inverse''
%map $k[V]\otimes V\to V$ sending $s\otimes u$ to an element that we
%denote by $s\backslash u$ and such that $s\backslash(sv)=v$. The
%torsion of a formal flat affine connection is a linear map
%$$T:k[V]\otimes V\otimes V\to V$$ defined as
%$$T(s)(y,z)=-[sy,sz].$$
%Its covariant derivatives $\nabla^n T(s)(x_1,\ldots,x_n,y,z)$ are
%linear maps $k[V]\otimes V^{\otimes n+2}\to V$ defined by the
%formula (\ref{eq:cov}) above.

Now, assume that a unital formal multiplication $F$ is given on $V$
and denote by $\mu_1\cdot \mu_2$ the corresponding product on
distributions. As mentioned in Section~\ref{subsec:mono}, it gives
rise to a formal connection simply by restricting $F$ to
$k[V]\otimes V$. The action of the adapted vector fields on
functions is easily derived from the product $\mu_1\cdot \mu_2$ on
$k[V]$.

\begin{lemma}
\label{lem:factorgamma} Let $\gamma\colon k[V] \to k[V]$ a linear
map that satisfies $$\Delta(\gamma(\mu)) = \sum\gamma(\mu_{(1)})
\otimes \mu_{(2)} + \mu_{(1)} \otimes \gamma(\mu_{(2)}).$$ Then
$$\gamma(\mu) = \sum  \mu_{(1)}
\pi_{V}\bigl(\gamma(\mu_{(2)})\bigr).$$
\end{lemma}
\begin{proof}
Let $S$ denote the antipode of the symmetric algebra $k[V]$
considered as a Hopf algebra. Since
\begin{align*}
\Delta\Bigl(\sum  S(\mu_{(1)}) \gamma(\mu_{(2)}) \Bigr) &= \sum
S(\mu_{(1)}) \gamma(\mu_{(2)})\otimes S(\mu_{(3)}) \mu_{(4)} \\
&\quad +
S(\mu_{(3)})
\mu_{(4)} \otimes S(\mu_{(1)}) \gamma(\mu_{(2)}) \\
&= \sum  S(\mu_{(1)}) \gamma(\mu_{(2)}) \otimes 1 + 1 \otimes
\sum  S(\mu_{(1)}) \gamma(\mu_{(2)})
\end{align*}
we have that $\sum  S(\mu_{(1)}) \gamma(\mu_{(2)}) $ is primitive.
Considering the degrees of the terms in this expression, we see that
$\sum S(\mu_{(1)}) \gamma(\mu_{(2)}) = \pi_V(\gamma(\mu))$. Thus
$$\gamma(\mu) = \sum \mu_{(1)}S(\mu_{(2)}) \gamma(\mu_{(3)})= \sum
\mu_{(1)} \pi_V(\gamma(\mu_{(2)})).$$

\end{proof}

\begin{lemma}
\label{lem:adapted_fields} For any $x\in V$ and $f \in k[V]^*$ we
have
$$
x^*(f)(\mu) = f(\mu\cdot x).
$$
\end{lemma}
\begin{proof}
Notice that $\gamma \colon \mu \mapsto \mu\cdot x$ satisfies
$\Delta(\gamma(\mu)) = \sum\gamma(\mu_{(1)}) \otimes \mu_{(2)} +
\mu_{(1)} \otimes \gamma(\mu_{(2)})$. This implies that
\begin{align*}
x^*(f)(\mu) = \sum  f(\mu_{(1)}x^*(\mu_{(2)})) = \sum  f(\mu_{(1)}\pi_V(\mu_{(2)}\cdot x)) = f(\mu\cdot x).
\end{align*}
\end{proof}

If $T$ is the torsion tensor of this connection, then setting
$$\langle x_1, \ldots, x_n; y,z \rangle_{F} =
\nabla_{x^*_1}\cdots \nabla_{x^*_n}T(y^*,z^*)(1)$$
we obtain an $n+2$-linear operation on $V$ for all $n\geq 0$. In
case that $G$ is an analytic local loop, the corresponding affine
flat connection is determined by its adapted vector fields $v^*$,
$v\in {\mathrm T}_e G$. For any analytic function $f$ on $G$ and any
distribution $\mu$ with support at the identity $e$, the
construction of Mikheev and Sabinin provides $\mu(v^*(f)) = \mu(g)$
with $g\colon a \mapsto v(f \circ L_a)$, so $\mu(v^*(f)) = (\mu
\cdot v) (f)$ in the bialgebra of distributions of $G$ with support
at the identity. Under the identification of analytic functions with
elements of $k[V]^*$, Lemma~\ref{lem:adapted_fields} shows that the
definition of adapted vector fields that we present agrees with this
one. Therefore, the formal connection, torsion and bracket
operations that we define agree with the corresponding constructions
by Mikheev and Sabinin. In \cite{SM87} they proved
\begin{proposition}
Assigning the set of operations $\langle x_1, \ldots, x_n, y,z
\rangle_{F}$ to a formal multiplication $F$ gives a functor from the
category of formal loops to that of Sabinin algebras with trivial
multioperator.
\end{proposition}

The torsion tensor also admits a simple interpretation in terms of
the product $\mu_1 \cdot \mu_2$.

\begin{lemma}
For any $x,y\in V$ and $\mu \in k[V]$ it holds
$$
T(x^*,y^*) (\mu) = \pi_V\left((\mu\cdot y)\cdot x - (\mu\cdot
x)\cdot y\right).
$$
\end{lemma}
\begin{proof}
We have that
\begin{align*}
T(x^*,y^*)(\mu) &= -[x^*,y^*](\mu) = \sum x^*(\mu_{(1)} y^*(\mu_{(2)})) - y^*(\mu_{(1)} x^*(\mu_{(2)})) \\
&= \sum  x^*(\mu_{(1)} \pi_V(\mu_{(2)}\cdot y))- y^*(\mu_{(1)} \pi_V(\mu_{(2)}\cdot x))
\\
&= \sum  x^*(\mu\cdot y) - y^*(\mu\cdot x) = \pi_V((\mu\cdot y)\cdot
x - (\mu\cdot x)\cdot y)
\end{align*}
\end{proof}

Recall that a set of multilinear brackets $\langle x_1, \ldots, x_n;
y,z \rangle$ on $V$ can be defined via the Shestakov-Umirbaev
operations.
\begin{theorem}\label{thm:brackets}
The operations $\langle x_1, \ldots, x_n; y,z \rangle$ of Shestakov
and Umirbaev identically coincide with  the operations $\langle x_1,
\ldots, x_n; y,z \rangle_F$ of Mikheev and Sabinin.
\end{theorem}
\begin{proof}
Let $\{v_i\}_i$ be a basis of $V$ and define formal functions $\{ f_i \}_i$ such that
$$
\langle \underline{\mu};y,z\rangle  = \sum_i f_i(\mu) v_i.
$$
The map $\gamma \colon \mu \mapsto (\mu \cdot z) \cdot y - (\mu
\cdot y) \cdot z$ satisfies the condition in
Lemma~\ref{lem:factorgamma} so
\begin{align*}
\sum  \mu_{(1)}T(y^*,z^*)(\mu_{(2)}) &= (\mu\cdot z) \cdot y -
(\mu\cdot y) \cdot z =
\sum  \mu_{(1)}\cdot \langle \underline{\mu}_{(2)}; y,z\rangle \\
&=\sum  \sum_i f_i(\mu_{(2)})\mu_{(1)} \cdot v_i=
\sum \sum_i  f_i(\mu_{(1)})\mu_{(2)}\pi_V(\mu_{(3)} \cdot v_i)\\
&= \sum \sum_i \mu_{(1)} (f_i(\mu_{(2)}) v_i^*(\mu_{(3)}))\\
&= \sum \mu_{(1)}\left(\sum_i f_iv^*_i(\mu_{(2)})\right).
\end{align*}
This proves that
$$
T(y^*,z^*) = \sum_i f_i v^*_i.
$$
The covariant differentiation of the torsion $T$ is then given by
$$
\nabla_{x^*_1} \cdots  \nabla_{x^*_n} T(y^*,z^*) = \sum_i
x^*_1\cdots x^*_n(f_i) v^*_i
$$
and the operations of Mikheev and Sabinin are recovered as
\begin{align*}
\langle x_1,\dots, x_n; y,z \rangle_F &= \nabla_{x^*_1} \cdots  \nabla_{x^*_n} T(y^*,z^*)(1)  \\
&=\sum_i f_i(((x_1\cdot x_2)\cdots )\cdot x_n) v_i \\
&= \langle x_1,\dots, x_n; y,z\rangle
\end{align*}
as desired.
\end{proof}

\iffalse

Let $F$ be unital formal multiplication on a vector space $V$.
Formal power series $k[V]\to V$ with no constant term can be
substituted as arguments into $F$ and this gives rise to a formal
loop $\widetilde{F}$ defined on the vector space of such formal
power series.  Given an arbitrary ordered $n$-tuple of formal series
$f_1,\ldots,f_n$ there is a linear homomorphism $f: {\mathcal
F}_n\to\widetilde{F}$ that sends the generator $x_i$ to $f_i$, and a
product of two non-associative monomials $m_1$ and $m_2$ to the
power series $F(f(m_1), f(m_2))$. Since  Theorem~\ref{thm:brackets}
holds for $ {\mathcal F}_n$, this implies that it also holds for
$\widetilde{F}$.

On the other hand, specifying a basis $\{e_{\alpha}\}$ of $V$ we get
an inclusion map of $V$ into $\mathrm{End}(V)$ which identifies $V$
with the subspace of diagonal matrices with finite number of
entries. In turn, $\mathrm{End}(V)$ is the same as {\em linear}
formal power series $k[V]\to V$, an we get an injective
``tautological'' homomorphism $F\to\widetilde{F}$. It embeds the
Sabinin algebra of $F$ into that of $\widetilde{F}$ and, hence
Theorem~\ref{thm:brackets} is true for $F$ as well.

\begin{remark}
In fact, it may be more natural to consider the formal loop
$\widetilde{F}$ of formal power series as a genuine loop, since the
product of two power series is again a power series. Similarly, the
formal loop ${\mathcal F}_n$ can be completed with respect to the
degree of the non-associative polynomials, and the resulting set of
power series also forms a loop.
\end{remark}
\fi

\subsection{Multioperators}
For a local analytic loop $(G,\cdot)$, the multioperator is a series
of operations on the tangent space $V=\mathrm{T}_eG$ at the identity
of $G$ given by
$$\Phi'(x,y)= \log{\Phi(\exp{x},\exp{y})}$$
where $\Phi$ is as in (\ref{eq:phi}). The homogeneous components of
$\Phi'$ are linear maps
$$\Phi'_{i,j}:k[V]_i\otimes k[V]_j\to V.$$
Each $\Phi'_{i,j}$ can be thought of either as a multilinear map
$V^{\otimes i+j}\to V$ which is totally symmetric in two groups of
variables, namely, the first $i$ and the last $j$ variables, or as a
polynomial map in two variables and bidegree $(i,j)$. In the
language of Section~\ref{subsec:mono} $\Phi'$ is a similarity.

The construction works for arbitrary formal loops if, instead of the
exponential map $\mathrm{T}_eG\to G$ one uses the exponential series
as defined in the Appendix. In particular, let us consider it for
the formal loop of {\em non-associative polynomials.}

Let $\mathcal{S}$ be a set and $k\{\mathcal{S}\}$ - the unital free
non-associative algebra generated by the elements of $\mathcal{S}$.
Denote by $\mR$ the ideal in $k\{\mathcal{S}\}$ generated by
$\mathcal{S}$. There is a unital formal multiplication on $\mR$
sending $x\otimes y$ to $x+y+xy,$ where the product $xy$ is taken in
$\mR$.

For $\alpha,\beta\in\mR$ write $a = \exp{\alpha}$, $b = \exp{\beta}$
and $\Phi'=\sum\Phi'_{i,j}(\alpha;\beta)$. (Here we treat
$\Phi'_{i,j}$ as a function of two variables $\alpha,\beta$ which is
of degree $i$ in $\alpha$ and $j$ in $\beta$.) Then (\ref{eq:phi})
has the form
$$
\exp_{a}{(a\Phi')}=\sum_{k = 0}^\infty {\frac{1}{k!}}
((a\Phi')\cdots)\Phi' = ab.
$$
This formula may be seen as a recursive definition of $\Phi'$. For
instance, expanding $a$ as a series in $\alpha$ we see that
\begin{align*}
 \Phi'_{1,3} &= -\frac{1}{12}[\beta,(\alpha,\beta,\beta)] -
 \frac{1}{6} p_{1,2}(\alpha;\beta^2;\beta)\\
 \Phi'_{2,3} &= -\frac{1}{12}(\alpha,\beta,(\alpha,\beta,\beta))
 + \frac{1}{12}(\alpha,(\alpha,\beta,\beta),\beta) \\
 &\qquad\qquad\qquad-\frac{1}{24}[v,p_{2,1}(\alpha^2;\beta;\beta)] -
 \frac{1}{12} p_{2,2}(\alpha^2;\beta^2;\beta).
\end{align*}
These expressions are essentially different from the multioperator of Shestakov and Umirbaev
$$\Phi_{i,j}^{SU}= p_{i,j-1}(\alpha^i,\beta^{j-1},\beta)/{i!j!}.$$
In general, we do not have such a closed formula for the
Sabinin-Mikheev multioperator.

%%%%%%%%%%%%%%%%%%%%%%%%%%%%%%%%%%%%%%%%%%%%%%%%%%%%%%%%%%%%%%%%%%%%%%%%%%%%%%%%
%%%%%%%%%%%%%%%%%%%%%%%%%%%%%%%%%%%%%%%%%%%%%%%%%%%%%%%%%%%%%%%%%%%%%%%%%%%%%%%%
%%%%%%%%%%%%%%%%                                                %%%%%%%%%%%%%%%%
%%%%%%%%%%%%%%%%                                                %%%%%%%%%%%%%%%%
%%%%%%%%%%%%%%%%              LINEAR FORMAL LOOPS               %%%%%%%%%%%%%%%%
%%%%%%%%%%%%%%%%                                                %%%%%%%%%%%%%%%%
%%%%%%%%%%%%%%%%                                                %%%%%%%%%%%%%%%%
%%%%%%%%%%%%%%%%%%%%%%%%%%%%%%%%%%%%%%%%%%%%%%%%%%%%%%%%%%%%%%%%%%%%%%%%%%%%%%%%
%%%%%%%%%%%%%%%%%%%%%%%%%%%%%%%%%%%%%%%%%%%%%%%%%%%%%%%%%%%%%%%%%%%%%%%%%%%%%%%%

\section{Linear formal loops}

Any finite--dimensional unital algebra $A$ over the real numbers
defines a local loop in a neighborhood on the identity $1$. By
translation $x\mapsto x-1$ we obtain a local loop  in a neighborhood
of $0$. The product $xy$ of this local loop is related with the
product $x*y$ of $A$ by
\[
xy = x + y + x*y.
\]
This formula, in fact, defines a unital formal multiplication on $A$
considered as a vector space. We shall denote this formal loop by
$G$. Note that the existence of the identity in $A$ is not relevant
here, so $A$ can be taken to be non-unital.
\iffalse When a basis $\{
e_1,\dots, e_n\}$ of $A$ is fixed we may identify $A$ with $k^n$.
The product $x*y$ induces a corresponding multiplication
\[
X*Y=(x_1,\dots,x_n)*(y_1,\dots,y_n) =(\sum_{i,j=1}^n\xi_{ij}^k x_i
y_j)_{k=1}^n
\]
with $\xi_{ij}^k \in k$, and a local loop on $e=(0,\dots,0)$
\[
G(X,Y) =XY= X+Y +X*Y.
\]
From a formal point of view $G(X,Y)$ is a formal multiplication and
it is natural to consider the relationship between the
Shestakov--Umirbaev operations on the Sabinin algebra $\Prim(k[G])$
and the corresponding operations on the algebra $(A,*)$.
\fi

As a vector space, $A$ can be identified with $\Prim(k[G])$ and,
hence, there are two ways to give the structure of a Sabinin algebra
to $A$: using the Shestakov-Umirbaev operations in $k[G]$ and in the
algebra $(A,*)$.

\begin{theorem}\label{thm:lin}
$\Prim(k[G])$ and $\SU(A,*)$ coincide as Sabinin algebras.
\end{theorem}
\begin{proof}
Let $(A^\#,*)= k1\oplus A$ the algebra obtained by adding a formal unit element $1$ to $A$.

The map $\pi_{A^\#}\colon  k[G]\to A^\#$ $\mu \mapsto \epsilon(\mu)1
+ \pi_A(\mu)$ which assigns to a distribution its component of
degree $\leq 1$ is, in fact, a homomorphism of algebras. Indeed,
$G(\mu_1,\mu_2) = \epsilon(\mu_2)\pi_A(\mu_1) +
\epsilon(\mu_1)\pi_A(\mu_2) + \pi_A(\mu_1)* \pi_A(\mu_2)$ so
$\pi_{A^\#}(G'(\mu_1,\mu_2)) =\epsilon(\mu_1)\epsilon(\mu_2) 1 +
G(\mu_1,\mu_2) = \epsilon(\mu_1)\epsilon(\mu_2) 1 +
\epsilon(\mu_2)\pi_A(\mu_1) + \epsilon(\mu_1)\pi_A(\mu_2) +
\pi_A(\mu_1)* \pi_A(\mu_2) = \pi_{A^\#}(\mu_1) * \pi_{A^\#}(\mu_2)$.
Since the Shestakov-Umirbaev operations are functorial with respect
to algebra homomorphisms then they coincide on $\Prim(k[G])=A$.

\end{proof}

\begin{corollary}
With the previous notation $k[G]\cong U(\SU(A))$ as bialgebras.
\end{corollary}

\begin{definition}
A formal loop $F$ is called { linear} if there exists a
finite-dimensional vector space $A$ with a bilinear product $x*y$
and a homomorphism
\[ \psi \colon F\rightarrow G\] with $G(\bx,\by) = \bx + \by + \bx*\by$,
where $\bx*\by$ stands for the formal map $k[A] \otimes k[A] \to A$,
$\mu_1\otimes \mu_2 \mapsto \pi_A(\mu_1)*\pi_A(\mu_2)$, such that
the induced $\Psi\colon k[G]^*\rightarrow k[F]^*$, $g\mapsto g \circ
\psi'$ is an epimorphism.
\end{definition}

\begin{lemma}
\label{lem:epi_mono} Let $F, G$ be  formal loops and $\psi \colon
F\rightarrow G$  a homomorphism of formal loops. Then $\Psi\colon
k[G]^*\rightarrow k[F]^*$, $g\mapsto g(\psi')$ is surjective if and
only if  $\psi'\colon k[F]\rightarrow k[G]$
is injective.
\end{lemma}

\begin{proposition}
Let $F$ be a formal loop. Then $F$ is linear if and only if there
exists a finite--codimensional ideal $I$ of the algebra $k[F]$ with
$I\cap \Prim(k[F]) = 0$.
\end{proposition}
\begin{proof}
Suppose that $F$ is linear. Then there exist $G(\bx,\by)= \bx+ \by+\bx*\by$, with
$x*y$ bilinear, and $\psi\colon F\rightarrow G$ homomorphism such
that the induced homomorphism $\psi'\colon k[F]\rightarrow k[G]$ is
injective. In the proof of Theorem~\ref{thm:lin} we saw that the
identity map on $A$ extends to surjective homomorphism $\pi_{A^\#}\colon
k[G]\rightarrow (A^\#, *)$. The kernel $I$ of the composition
$\pi_{A^\#}\psi'$ is a finite--codimensional ideal of $k[F]$ with $I\cap
\Prim(k[F]) =0$.

Conversely, assume that $I$ is a finite--codimensional ideal of
$k[F]$ such that $I\cap \Prim(k[F])=0$. Consider $A=k[F]/I$, with
the product denoted by $*$, and $G(\bx,\by) = \bx + \by + \bx*\by$  the
corresponding unital formal multiplication. Using the universal
property of $k[F]$ we see that the monomorphism $\Prim(k[F])
\rightarrow A\cong \Prim(k[G])$ of Sabinin algebras induces a
homomorphism of bialgebras $\psi'\colon k[F]\rightarrow k[G]$ with
injective restriction to $\Prim(k[F])$. By Theorem~2.4.11 in
\cite{Ab80}, this map $\psi'$ must be injective too. The proposition
now follows from Lemma~\ref{lem:epi_mono}.
\end{proof}

In \cite{PS04} it was proved that any Moufang formal loop is linear,
a result that extends Ado's theorem to Malcev algebras. However,
there exist formal multiplications that are not linear. Important
examples come from Bruck loops. A {\em Bruck loop} is a loop that
satisfies the Bol identity
\[
a(b(ac))=(a(ba))c
\]
(which implies that $L^{-1}_a=L_{a^{-1}}$  for some $a^{-1}$)  and
the \emph{automorphic inverse property}
\[
(ab)^{-1}=a^{-1}b^{-1}
\]
for all $a,b,c$. The bialgebra of distributions $k[F]$ of a formal
Bruck loop satisfies
\[
\sum \mu_{(1)}(\nu(\mu_{(2)}\eta))=\sum
(\mu_{(1)}(\nu\mu_{(2)}))\eta
\]
and there exists a map $S$ such that
\[
\mu\backslash \nu =S(\mu)\nu
\]
and
\[
S(\mu\nu)=S(\mu)S(\nu)
\]
for all $\mu,\nu,\eta\in k[F]$. In this case all the operations of
the Sabinin algebra $\Prim(k[F])$ are easily described in terms of a
Lie triple system, and conversely, any Lie triple system provides a
formal Bruck loop law. However, Lie triple systems that are not
nilpotent do not provide linear multiplications \cite{Pe08}.

%%%%%%%%%%
\subsection{Some examples}
%%%%%%%%%%

Let us consider the bilinear product on $k^3$
\[
x*y =(x_1 y_1 + x_2 y_3 + x_3 y_2, x_1 y_2 + x_2 y_1, x_1 y_3 + x_3
y_1)
\]
and the formal multiplications
\begin{eqnarray*}
G(x,y) &=& x+y+x*y\\
F((x_2,x_3),(y_2,y_3))&=&\frac{1}{1+x_2 y_3 + x_3 y_2}(x_2 + y_2,
x_3+y_3)
\end{eqnarray*}
The map
\[
\phi=\Bigl(\frac{x_2}{1+x_1},\frac{x_3}{1+x_1}\Bigr)
\]
defines a homomorphism $\phi\colon G \rightarrow F$ of formal loops.
It induces a surjective homomorphism $\phi'\colon k[G]\rightarrow
k[F]$ determined by $\phi'(\partial_1) =0$, $\phi'(\partial_2)
=\partial_2$ and $\phi'(\partial_3) = \partial_3$, where
$\partial_i$ is the basis vector of $k^3$ corresponding to the
coordinate $x_i$. These formulae come from considering $(V,(\,,\,))$
a two--dimensional vector space with a bilinear form of maximal Witt
index, $A = \mathbb{R} e\oplus V$ the Jordan algebra with the
product
\[
(\alpha e + a)(\beta e + b) = (\alpha\beta +(a,b))e + \alpha b +
\beta a
\]
and the formal loop determined by $A$. The subspace $\mathbb{R}e$
may be thought of as a normal subloop  and $F$ as the quotient of
$G$ by $\mathbb{R}e$. Although the formal loop $G$ is linear, we
shall see that $F$ is not . To simplify the notation involved in our
computations, we shall identify $k[G]$ with $U(\SU(A))$.

\begin{lemma}
\label{lem:Jordan} Let $A$ be a Jordan algebra, then in $U(\SU(A))$
we have
\[
p(a,\underline{xc},b) = -\sum p(c,\underline{x}
_{(1)},p(a,\underline{x} _{(2)},b)) + \epsilon(x)(a,c,b)
\]
and
\[
(a,xc,b) = (a,x,b)c+x(a,c,b)
\]
for any primitive $a,b,c$ and any $x\in U(\SU(A))$.
\end{lemma}
\begin{proof}
The map $x\mapsto (a,x,b)$ is a derivation of any Jordan algebra so
in $A$
\[
(a,xc,b)= \left\{ \begin{array}{ll} \sum x_{(1)}p(a,\underline{x}
_{(2)}c,b) +
 \sum (x_{(1)}c)
p(a,\underline{x} _{(2)},b)\\
\\
(a,x,b)c + x(a,c,b) =  \sum c\big(x_{(1)}p(a,\underline{x}
_{(2)},b)\big) + x(a,c,b)
\end{array}\right.
\]
thus
\begin{eqnarray*}
 \sum x_{(1)}p(a,\underline{x}_{(2)}c,b) \hskip -2pt&=&  \sum -(x_{(1)},c,p(a,\underline{x} _{(2)},b))
 + x(a,c,b)\\
&=&  \sum
-x_{(1)}p(\underline{x}_{(2)},c,p(a,\underline{x}_{(3)},b))  +
 \sum x_{(1)}\epsilon(x_{(2)}) (a,c,b)
\end{eqnarray*}
Dividing by $x_{(1)}$ we get the first equality. The second equality
follows from the first one by reversing our argument in $U(\SU(A))$
(notice that $U(\SU(A))$ is commutative).
\end{proof}

\begin{theorem}
The formal loop \[ F((x_2,x_3),(y_2,y_3))=\frac{1}{1+x_2 y_3 + x_3
y_2}(x_2 + y_2, x_3+y_3)
\]
is not linear.
\end{theorem}
\begin{proof}
Any finite--codimensional ideal of $k[F]$ that meets trivially the
primitive elements provides a finite--codimensional ideal of $k[G]$
that contains $\partial_1$ and with trivial intersection with
$k\partial_2 +k\partial_3$. With the identification $k[G]\cong
U(\SU(A))$ we obtain a finite--codimensional ideal $I$ of
$U=U(\SU(A))$ with $e\in I$ and $V\cap I = 0$. Let us show that this
is not possible. We will fix $a,b\in V$ with $(a,a) = 0=(b,b)$ and
$(a,b) = 2$.

Since $A$ is a Jordan algebra, $A$ is commutative and
power--associative, so the formal loop determined by $A$ also is.
The universal enveloping algebra $U$ is commutative and the powers
$x^n$ are well--defined for any $x\in \Prim(U)$. The dimension of
$U/I$ is finite so we can find a linear combination $a^N +\alpha_1
a^{N-1} +\cdots + \alpha_{N-1} a \in I$. By the previous Lemma, we
conclude that $a^N\in I$. We also assume that $N$ is minimal with
respect to this property.

In $A$ the powers $a^m$ vanish if $m\geq 2$. In such a case the
relation $ (a^m,b,b)=\sum a^m_{(1)}p(\underline{a}^m_{(2)},b,b)$
implies that
\[
p(\underline{a}^m,b,b) = -map(\underline{a}^{m-1},b,b) = \cdots =
(-1)^{m-1}4 m! a^{m-2}
\]
and we obtain $p(a,b,b) = 2b, p(\underline{a}^2,b,b) = -8 e,
p(\underline{a}^3,b,b)= 24 a$ and $p(\underline{a}^m,b,b) =0$ if
$m\geq 4$. Let us use these formulae to prove that $a^N\in I$
implies $a\in I$, which is not possible because $V\cap I = 0$. In
case that $N=2$, in $U(A)$ we have that modulo $I$
\begin{eqnarray*}
0&\equiv &(a^2,b,b)a = p(\underline{a}^2,b,b) a + 2ap(a,b,b)\cdot a
= -8ea + 4ab\cdot a \\ &\equiv& 4a^2b -4(a,a,b)\equiv -4(a,a,b) =
8a.
\end{eqnarray*}
In case that $N\geq 3$ then
\begin{eqnarray*}
0&\equiv & (a^N,b,b) a = N a^{N-1}p(a,b,b)\cdot a + \binom{N}{2}
a^{N-2}p(\underline{a}^2,b,b)\cdot a \\
&& + \binom{N}{3} a^{N-3}p(\underline{a}^3,b,b)\cdot a  \\
&=& 2N a^{N-1}b\cdot a -4 N(N-1) a^{N-2}e \cdot a + 4N(N-1)(N-2)
a^{N-1}\\
&\equiv& -2N (a,a^{N-1},b) + 4N(N-1)(N-2) a^{N-1} \\
&=& 4N(N-1)^2a^{N-1}
\end{eqnarray*}
so $a^{N-1}\in I$, a contradiction with the minimality of $N$.
\end{proof}

Operations $\langle\,;\,,\,\rangle$ on Jordan algebras are
determined by a Lie triple system. The same relation holds for Bol
algebras with trivial binary product. This indicates that a formal
loop determined by a Jordan algebra is similar to a formal Bruck
loop \cite{Pe07}.

\begin{proposition}
If $A$ is a Jordan algebra, in $\SU(A)$ we have
\[
\langle \underline{xc}; a,b\rangle =  \sum \langle
\underline{x}_{(1)}; c,\langle \underline{x}_{(2)}; a,b\rangle
\rangle \quad\text{and}\quad \langle c;a,b\rangle =-(a,c,b)\] if
$\vert x\vert \geq 1$.
\end{proposition}
\begin{proof}
Since $U(\SU(A))$ is commutative, by Lemma \ref{lem:Jordan}
\begin{eqnarray*}
(xc,a,b) -(xc,b,a) &=& (xc)a\cdot b -(xc)(ab) -(xc)b\cdot a
+(xc)(ba)\\ &=& -(b,xc,a) = -(b,x,a)c-x(b,c,a)
\end{eqnarray*}
and
\begin{eqnarray*}
 \sum (xc)_{(1)}\cdot \langle \underline{xc}_{(2)};b,a\rangle &=& (xc,a,b) -(xc,b,a)
= -(b,x,a)c-x(b,c,a) \\ &=&  \sum  x_{(1)}\langle
\underline{x}_{(2)};b,a\rangle\cdot c -x(b,c,a)
\end{eqnarray*}
so
\begin{eqnarray*}
 \sum x_{(1)}\cdot \langle \underline{x}_{(2)}c;b,a\rangle &=&   \sum -x_{(1)}c\cdot \langle
\underline{x}_{(2)};b,a\rangle  \\
&& +  \sum x_{(1)}\langle
\underline{x}_{(2)};b,a\rangle\cdot c -x(b,c,a)\\
&=&  \sum -(c,x_{(1)},\langle \underline{x}_{(2)};b,a\rangle) -x(b,c,a)\\
&=&  \sum x_{(1)}\langle \underline{x}_{(2)}; c,\langle
\underline{x}_{(3)};b,a\rangle\rangle-x(b,c,a)
\end{eqnarray*}
as desired.
\end{proof}

\appendix
\section*{Non-associative exponential and logarithm}
\subsection*{The exponential}

Let $\mRhat$ be the algebra of non-associative power series in some
set of variables with coefficients in  $k$  and with no constant
term. Given $X\in \mRhat$ we define $\exp{X}\in 1+\mRhat$ as
\[\exp{X}=1+X+\frac{X^2}{2!}+\frac{X^2X}{3!}+\frac{((X^2)X)X}{4!}+\ldots\]
It is readily seen that $\exp{X}$ is the value at $t=1$ of the
solution of the differential equation
\[\frac{da}{dt}=aX\]
with the initial condition $a(0)=1$.

One may think of the algebra $\mRhat$ is the tangent space at 1 to
the multiplicative loop $1+\mRhat$. Right multiplication by $b\in
1+\mRhat$ defines a parallel transport of $\mRhat$ to $b+\mRhat$.
More generally, the canonical connection on $1+\mRhat$ is defined by
transporting $b+X\in b+\mRhat$ to $c+c(b\backslash X)\in c+\mRhat$
for all $b,c\in 1+\mRhat$.

Curves of the form $\exp{Xt}$ are the geodesics of the canonical
connection that pass through 1. It is equally easy to write down the
geodesics that pass through $b\in 1+\mRhat$. For $X\in \mRhat$
define $\exp_b{X}$ as
\[\exp_b{X}=b+X+\frac{X(b\backslash X)}{2!}+\frac{(X(b\backslash X))(b\backslash X)}{3!}+
\frac{((X(b\backslash X))(b\backslash X))(b\backslash
X)}{4!}+\ldots.\] Then $\exp_b{Xt}$ satisfies the equation
\[\frac{da}{dt}=a(b\backslash X)\]
with the initial condition $a(0)=b$.

It is easily verified that, just as in the associative case, $X\in
\mRhat$ is primitive if and only if $\exp{X}\in 1+\mRhat$ is
group-like, that is, $$\Delta \exp{X}=\exp{X}\otimes\exp{X}.$$ This
property, however, does not define the exponential series uniquely;
see, for instance, \cite{GH}.

\subsection*{The logarithm}

The power series $\log(1+x)$ is defined by $\exp(\log(1+x))=1+x$.
The coefficients of $\log(1+x)$ can be found as follows.

Assume that $\mRhat=\mRhat(x)$, the algebra of non-associative power
series in one variable $x$. (One can forget altogether about the
variable and think of the non-associative monomials in $x$ as of
rooted binary plane trees.)

Write $X=\sum_{\tau} X_{\tau} \tau$ where the sum runs over all
non-associative monomials $\tau$. Then $\exp X$ can be written as
\[\exp{X}=\sum_{\tau=(\ldots(\tau_1\tau_2)\ldots)\tau_k}
\frac{X_{\tau_1}X_{\tau_2}\ldots X_{\tau_k}}{k!}\cdot\tau.\] Writing
$\exp{X}=\sum_{\tau} a_{\tau} \tau$ we have
\begin{align*}
a_{(..(\tau_1\tau_2)\ldots)\tau_k}&=
X_{(..(\tau_1\tau_2)\ldots)\tau_k}+
\frac{1}{2!}X_{(..(\tau_1\tau_2)\ldots)\tau_{k-1}}X_{\tau_k}
\\&\quad+ \frac{1}{3!}X_{(..(\tau_1\tau_2
)\ldots)\tau_{k-2}}X_{\tau_{k-1}}X_{\tau_k}+\ldots.
\end{align*}
Also,
\begin{align*}
a_{(..(\tau_1\tau_2)\ldots)\tau_{k-1}}X_{\tau_k}&=
X_{(..(\tau_1\tau_2)\ldots)\tau_{k-1}}X_{\tau_k}+
\frac{1}{2!}X_{(..(\tau_1\tau_2)\ldots)\tau_{k-2}}X_{\tau_{k-1}}X_{\tau_k}+
\ldots,\\
a_{(..(\tau_1\tau_2)\ldots)\tau_{k-2}}X_{\tau_{k-1}}X_{\tau_k}&=X_{(..(\tau_1\tau_2)\ldots)\tau_{k-2}}X_{\tau_{k-1}}X_{\tau_k}\\
&\qquad\qquad+\frac{1}{2!}X_{(..(\tau_1\tau_2)\ldots)\tau_{k-3}}X_{\tau_{k-2}}X
_{\tau_{k-1}}X_{\tau_k}+ \ldots
\end{align*}
and so on.

Recall that the Bernoulli numbers $B_k$ satisfy the identity
\[\sum_{k=0}^{n-1}\frac{B_k}{k!(n-k)!}=0.\]
It follows that
\begin{multline*}
X_{(..(\tau_1\tau_2)\ldots)\tau_k}=
a_{(..(\tau_1\tau_2)\ldots)\tau_k}+
\frac{B_1}{1!}a_{(..(\tau_1\tau_2)\ldots)\tau_{k-1}}X_{\tau_k}\\+
\frac{B_2}{2!}a_{(..(\tau_1\tau_2)\ldots)\tau_{k-2}}X_{\tau_{k-1}}X_{\tau_k}+\ldots+
\frac{B_{k-1}}{(k-1)!}a_{\tau_1}X_{\tau_2}\ldots X_{\tau_k}.
\end{multline*}

Now, set $a_x=1$ and $a_{\tau}=0$ for $\tau\neq x$. Then the
$X_{\tau}$ are the coefficients of the power series $\log{(1+x)}$.
Setting $\tau=(..((x\tau_1)\tau_2)\ldots)\tau_k$ we see that
\[X_{\tau}=\frac{B_{k}}{k!}X_{\tau_1}\ldots X_{\tau_k}.\]

Given a binary rooted plane tree $\tau$ define $B_{\tau}$ and
$\tau!$ inductively as follows.

For $\tau=x$ set $B_\tau=\tau!=1$. If $\tau\neq x$, there is only
one way of writing $\tau$ as a product $(\ldots((x
\tau_1)\tau_2)\ldots)\tau_k$. Set
\[B_{\tau}=B_k\cdot B_{\tau_1}\ldots B_{\tau_k}\]
and
\[\tau! = k! \tau_1 !\ldots \tau_k !.\]
With this notation we have

\[\log(1+x)=\sum_{\tau} \frac{B_\tau}{\tau!} \cdot \tau. \]

\subsection*{Identities related to sums over trees}

This expression for the coefficients of the non-associative
logarithm implies certain identities on Bernoulli numbers. Imposing
the associativity condition on $\mRhat$, we turn our exponential
into the usual exponential series; therefore, our logarithm becomes
the usual logarithm. All monomials $\tau$ with $\deg \tau=n$ are
sent to the monomial $x^n$. We obtain
\begin{equation}\label{eq:bernoulli}
\sum_{\deg\tau=n} \frac{B_\tau}{\tau!}= \frac{(-1)^{n+1}}{n}.
\end{equation}

\def\b{\beta}
\def\l{\lambda}
\def\t{\tau}
\def\={\;=\;}
\def\B{\mathcal B}

A direct proof of (\ref{eq:bernoulli}), together with a
generalization of it, was communicated to us by D.\ Zagier.

Choose arbitrary weights $\b_1,\,\b_2,\dots$  and for $n\geq 1$
define $\l_n$ as $\sum\b_\t$, where the sum runs over plane rooted
trees $\t$ of degree~$n$ and $\b_\t$ is defined as
$\b_{i_1}\cdots\b_{i_k}$ if the vertices of $\t$ have
$i_1,\dots,i_k$ outgoing branches. Since each such tree consists of
a root which is joined to the roots of some (ordered) collection of
plane rooted trees, say $\t_1,\dots,\t_r$ of degrees
$n_1,\dots,n_r\ge1$ with $\sum_in_i=n-1$, we have
$$\l_1=1\,,\qquad
\l_n\=\sum_{r\ge1}\;\,\b_r\!\!\!\sum_{\stackrel{n_1,\dots,n_r\ge1}{n_1+\cdots+n_r=n-1}}\lambda_{n_1}\cdots\l_{n_r}\qquad\text{if
$n>1\,$.}$$ Hence the generating function
$L=L(x)=\sum\limits_{n=1}^\infty\l_nx^n$ satisfies the functional
equation
$$L\=x\,\biggl(1+\sum_{r=1}^\infty\b_r\,L^r\biggr),$$ or
$$\frac{L}{\B(L)}\=x\,,$$ where
$\B(t)=1+\sum\limits_{r=1}^\infty\b_r\,t^r\,$. For instance, if all
$\b_r=1$ then $\B(t)=\dfrac1{1-t}$ and therefore $L(1-L)=x$ or
$L=\frac12(1-\sqrt{1-4x})$, the standard generating function for the
number $\binom{2n}n/(n+1)$ of plane rooted trees of degree $n$ (=
number of length~$n$ bracketings = $n$th Catalan number).  If
$\b_r=B_r/r!$ then we have instead $\B(t)=\dfrac t{e^t-1}$ and hence
$x=e^{L}-1$ or $L=\log(1+x)$, giving  $\l_n=(-1)^{n-1}/n$, that is,
the formula (\ref{eq:bernoulli}).

\subsection*{Acknowledgments}
We are grateful to D.\ Zagier for his comments on the
formula~(\ref{eq:bernoulli}). The first author would like to thank Max-Planck-Institut f\"ur
Mathematik, Bonn for hospitality. Both authors were supported by the
SEP-CONACyT grant 44100, and the second author acknowledges support of the grants
MTM2007-67884-C04-03 and ANGI2005/05.

\bibliography{Biblio}
\bibliographystyle{amsalpha}

\end{document}